\documentclass[11pt,reqno]{amsart}

\usepackage{pgfplots}
\usepgfplotslibrary{colormaps,fillbetween}
\usepackage{mathrsfs}
\usetikzlibrary{arrows}
\usepackage{a4wide}
\usepackage{changepage}
\usepackage{amsmath,amssymb} 
\usepackage{bbm} 
\usepackage{color}
\usepackage{tikz}
\usepackage{tkz-graph}
\usepackage{tikz-cd} 
\usetikzlibrary{decorations.pathreplacing,calligraphy}
\usepackage{dsfont}
\usepackage{amscd}
\usepackage{graphicx}
\usepackage{mathtools}
\usepackage{todonotes}
\usepackage{enumerate}
\usepackage[normalem]{ulem}
\usepackage{float}
\usepackage{url}

\usepackage{scalerel,stackengine}

\usepackage{comment}

\usepackage[pagewise]{lineno}

\usepackage{bm}

\theoremstyle{plain}
\newtheorem{theorem}{Theorem}[section]
\newtheorem{prop}[theorem]{Proposition}
\newtheorem{lemma}[theorem]{Lemma}

\theoremstyle{definition}

\newtheorem{remark}[theorem]{Remark}

\newcommand{\R}{{\mathbb R}}
\newcommand{\N}{{\mathbb N}}

\newcommand{\X}{{\mathbb X}}
\newcommand{\T}{{\mathbb T}}

\newcommand{\mc}{\mathcal}

\renewcommand{\emptyset}{\varnothing}

\newcommand{\TT}{\mathbb{T}}

\usepackage{mathtools}

\stackMath
\newcommand\reallywidecheck[1]{%
\savestack{\tmpbox}{\stretchto{%
   \scaleto{%
\scalerel*[\widthof{\ensuremath{#1}}]{\kern-.6pt\bigwedge\kern-.6pt}%
     {\rule[-\textheight/2]{1ex}{\textheight}}
   }{\textheight}%
}{0.5ex}}%
\stackon[1pt]{#1}{\scalebox{-1}{\tmpbox}}%
}

\begin{document}

\title[Large deviations and differentiability of pressure]{Large deviations for observables with logarithmic singularity and differentiability of the pressure function
}

\author{Matthew Nicol}
\address{Matthew Nicol\\ Department of Mathematics\\
University of Houston\\
Houston\\
TX 77204\\
USA} \email{nicol@math.uh.edu}
\urladdr{\url{http://www.math.uh.edu/~nicol/}}

 \author{Tanja I.\ Schindler}
 \address{Faculty of Mathematics and Computer Science, Jagiellonian University\newline 
 \hspace*{\parindent} Ul.\ Łojasiewicza 6, 30-348 Krakow, Poland\newline 
 Department of Mathematics and Statistics,
 University of Exeter\newline 
 \hspace*{\parindent} Harrison Building, Exeter EX4 4QF,
 United Kingdom}
 \email{tanja.schindler@uj.edu.pl, t.schindler@exeter.ac.uk}
 \urladdr{http://tanja-schindler.github.io}

 \makeatother
 
 \thanks{
 TS was partly supported by the Austrian Science
Fund FWF: P 33943-N,  by a grant from the priority research area SciMat under the
Strategic Programme Excellence Initiative at Jagiellonian University
and by the European
Union’s Horizon Europe research and innovation programme under the Marie Sklodowska-Curie grant agreement ErgodicHyperbolic - 10115118.
 MN  would like to thank  the Leverhulme Trust  for support and the University of Warwick for hospitality. MN and TS would like to thank
the Schr\"odinger Institute for support through the 2024 Workshop on Extremes and Rare Events.  
Moreover, the authors would like to thank the Universities of Exeter, Houston and Warwick for their hospitality where parts of the discussions have taken place. 
Finally, the authors would like to thank Philipp Gohlke for useful discussions.
}



\begin{abstract}

We show that  logarithmic observables of form $\phi(x)=\log(|x-b|)$ on
the doubling map $T$
have exponential large deviations with
respect to the invariant measure if $b$ is preperiodic but not periodic with respect to $T$ giving a counterpart to the known result of strictly stretched exponential large deviations in the case that $b$ is
periodic.  These results are established by showing that the topological pressure function developed for unbounded observables is differentiable at $0$. 
\end{abstract}

\maketitle

\section{Introduction and notation}
Let $(T,X,\mu)$ be an ergodic probability-preserving dynamical system. If $\phi$ is an integrable observable on $X$ we say 
$S_n (\phi):=\sum_{j=0}^{n-1} \phi \circ T^j$ has exponential large deviations if for each $\varepsilon>0$ there exists
$C(\varepsilon)>0$ and $g(\varepsilon)>0$ such that for all $n\geqslant 0$
\[
\mu \Big(\Big|\frac{1}{n}S_n (\varphi)-\overline\varphi\Big| > \varepsilon\Big)\leqslant C(\varepsilon)  e^{-ng(\varepsilon)}.
\]

If for each $\varepsilon>0$ there exist
$C(\varepsilon) >0$ and $0<\gamma (\varepsilon)  <1$ such that for all $n\geqslant 0$
\[
\mu \Big(\Big|\frac{1}{n}S_n (\phi)- \mu (\phi) \Big| > \varepsilon\Big)\leqslant C(\varepsilon)  e^{-n^{\gamma(\varepsilon)}}
\]
we say $S_n (\phi)$ has stretched exponential large deviations of order $\gamma>0$.

We now recall the definition of rate function and some other notions of classical  large deviations theory pertaining to our results.   
A mean-zero integrable function $\phi: X \to \R$ on a dynamical system $(T,X,\mu)$
  is said to satisfy a large deviation principle with rate function
  $I(\alpha)$, if there exists a non-empty neighborhood $U$ of $0$ and a
  strictly convex function $I:U\to \R$, non-negative and vanishing
  only at $\alpha=0$, such that
\begin{align}\label{rate+}
\lim_{n\to\infty} \frac 1n\log \left(\mu (S_n (\phi)  \geqslant n \alpha)\right)& = -I(\alpha)
\end{align}
for all $\alpha>0$ in $U$ and
\begin{align}\label{rate-}
 \lim_{n\to\infty} \frac 1n\log \left(\mu (S_n (\phi)  \leqslant n \alpha) \right)&= -I(\alpha)
 \end{align}
for all $\alpha<0$ in $U$.

In the literature this is referred to as a first level or local (near the average) large deviations principle. There is a well-developed theory establishing exponential large deviations with a rate function 
for sufficiently regular (usually H\"older) observables on uniformly hyperbolic dynamical systems going back to Kifer and Young, \cite{kifer,You90}. 
 If $\phi$ is a mean-zero continuous observable on an SRB attractor then
 $\phi$ has exponential large deviations (see~\cite[Theorem 2
 (2)]{You90}).
 For mean-zero H\"older observables on Young Towers with exponential tails
 (we refer to~\cite{You90} or~\cite[Sections 2 and 4]{MN08} for the
 definition)
 which are not $L^1$ coboundaries in the sense that
 $\phi\not = \psi\circ T-\psi$ for any $\psi\in L^1 (\mu)$ such an
 exponential large deviations result holds with rate function
 $I_{\phi} (\alpha)$~\cite{MN08,RBY08}. A formula for the width
 of $U$ is given in \cite{RBY08} following a standard approach
 but it is not useful in concrete estimates.

 In general little is known about exponential large deviation estimates  for unbounded observables, even in the uniformly hyperbolic case. 
 A special case is considered in \cite{takahasi_cfs}  where the observable is piecewise constant on the partition generated by the underlying transformation on the interval $T$. 
 A morally similar situation occurs when studying countable Markov shifts with finite  variation on each cylinder set, see e.g.\  \cite{iommi, takahasi_markov}.

  For certain  non-uniformly expanding maps Ara\'ujo~\cite[Theorem
  3.1]{Ara07} gives an asymptotic condition, called exponential slow recurrence to $p$, on observables of form
  $\phi (x)=-\log (|x-p|)$
  which implies exponential large deviations for $\phi$. He also proves
  that singular points for certain Lorenz-like maps with unbounded
  derivatives satisfy this condition, which is in general difficult to
  check.
  These results were later also generalized to Teichm\"uller flows.

However, in general, the picture is far from being complete when the observable is unbounded.
In general, one should not 
expect exponential large deviations, even for uniformly expanding systems. It was shown in~\cite{NT_LD} that only stretched exponential decay of correlations holds for certain observables with a logarithmic singularity over the piecewise expanding maps. 

To formalize the situation and adapt it to our setting, let $T:[0,1]\to[0,1]$ be the doubling map given by $T(x)=(2x)$ (mod $1$).

 Moreover, let $\phi^b : [0,1] \to \R$ be an observable with a logarithmic singularity, of form
\begin{align*}
 \phi^b (x)= 2\log (|x-b|) 
\end{align*}
where $b\in [0,1]$.  We normalize and write $\overline{\phi}^b=\phi^b-\mu(\phi^b)$.
The factor $2$ is a normalisation which matches $\phi^b$ to the potential $\psi^b$ introduced in Section \ref{sec: diff pressure} that originates from the study of generalized Thue--Morse diffraction measures (see Remark~\ref{rem: Thue Morse pot}) and whose pressure function we study later.
We define the ergodic sum 
$S_n(\overline{\phi}^b) :=\sum_{k=0}^{n-1} \overline{\phi}^b \circ T^k$.  Given $\varepsilon>0$ we are interested in the rate of decay of 
\[
\mu (|S_n(\overline{\phi}^b)|\geqslant n\varepsilon).
\]
Results of~\cite{NT_LD} show  that if $b$ is periodic with respect to $T$ then the decay is at most stretched exponential. (In fact in~\cite{NT_LD} they consider $-\log (|x-b|)$ 
but their results  hold with a straightforward modification for $C \log (|x-b|)$ for any $C>0$.)
In \cite{pontes} it was shown that for any $b$ the decay is at least stretched exponential, but no lower bound was given. In recent work~\cite{defaveri_smith} 
 large deviations  for the observable $-\log (x)$  on  the doubling map were examined, and  detailed but non-rigorous predictions for the  distribution  were given as well as a transition at the scaling $n^{2/3}$.

In this paper we  show that for all  preperiodic but not periodic $b$ the decay rate is
exponential. This answers a question poised in~\cite{NT_LD}.  

One of the main tools for proving this result is to show that the topological pressure function defined in \cite{GKS} for an observable with logarithmic singularity is differentiable at $0$. As in this case the pressure function is no longer subadditive, the differentiability is a non-trivial problem; for similar differentiability results in other settings see Remark~\ref{rem: Thue Morse pot}.

The organization of this paper is as follows. In Section~\ref{results} we state our main results. This includes an improvement of the 
large deviations bound in the case that $b$ is periodic given in~ \cite[Theorem 3.4]{NT_LD}.
In Sections~\ref{sec: sketch} to \ref{sec: proofs large dev} we give the proof of the main theorem and in Section~\ref{sec: periodic case} we give the proof of the case $b$ being periodic.

\section{Statement of main results}\label{results}

\subsection{Observables with a logarithmic singularity at a non-periodic preperiodic point}
We consider $([0,1],\mathcal{B}, \mu, T)$ with $T(x)=2x\mod 1$, $\mu$ the Lebesgue measure, 
and $\phi^b(x):=2\log (|x-b|)$, for different values of $b\in [0,1)$ and set $\overline{\phi}^b:=\phi^b-\mu(\phi^b)$. 

\begin{theorem}\label{LDP}
 Let $b$ be preperiodic but not periodic, then $\overline{\phi}^b$  has exponential decay of  deviations. More precisely, there exists an interval $U=(0,u)$, $u>0$ and functions
 $C, I:U\to \mathbb{R}_{>0}$$C(\varepsilon)>0$ and function  $I>0$
such that if $\varepsilon\in U$, we have
\begin{equation}
 \mu(|S_n(\overline{\phi}^b) |\geqslant n\varepsilon) <C(\varepsilon)  e^{-n I  (\varepsilon)},\label{eq: main ineq}
\end{equation}
for all $n\geqslant 0$.
\end{theorem}

The main idea of the proof of this theorem will be given in Section \ref{sec: sketch}
with its lemmas proven in Section \ref{sec: diff pressure} and \ref{sec: proofs large dev}.

\subsection{Observables with a logarithmic singularity at a periodic point}\label{periodic}

We obtain an improvement to \cite[Theorem 3.4]{NT_LD}. To keep in line with that paper we consider 
functions with logarithmic singularities bounded below. Note that for $\alpha=1$ we have $\varphi_\mathsf{p}=-\phi^\mathsf{p}/2$, so that the upper tail estimated here corresponds to the lower tail of $\overline{\phi}^b$ studied in the rest of the paper.

\begin{prop}\label{prop: periodic part}
 Let $([0,1],\mathcal{B}, \mu, T)$ be a topologically mixing Rychlik map of the unit interval with $\mu$ the unique ergodic invariant probability measure, absolutely continuous with respect to the Lebesgue measure. Moreover, let $\varphi_\mathsf{p}(x):=  (-\log (|x-\mathsf{p}|))^{\alpha}$, 
 where $\mathsf{p}$ is periodic for $T$ and $\alpha > 0$. Assume that $T$ is $\mathcal{C}^1$ on a neighborhood of the orbit of $\mathsf{p}$. 
 Then for all sufficiently small $\varepsilon>0$ there exist $c(\varepsilon), C(\varepsilon), r(\varepsilon), R(\varepsilon) > 0$  and $n_0 (\varepsilon)$ such that for all $n \geqslant n_0 (\varepsilon)$  we have
\begin{align*}
 c(\varepsilon) \exp\left(-r(\varepsilon)n^{\frac{1}{1+\alpha}}\right)
 \leqslant \mu\left(S_n(\varphi_\mathsf{p}) - n\mu(\varphi_\mathsf{p}) \geqslant n\varepsilon\right)
 \leqslant C(\varepsilon)  \exp\left(-R(\varepsilon) n^{\frac{1}{1+\alpha}}\right).
\end{align*}
\end{prop}
We postpone the short proof to the end of the paper in Section \ref{sec: periodic case}.

\section{Skeleton of proof of Theorem~\ref{LDP}}\label{sec: sketch}

In this Section we will give the main steps to prove Theorem~\ref{LDP} by formulating some propositions and lemmas whose proof will be postponed to later sections.

We will first prove the results for the function 
\begin{align*}
 \psi^b(x)&:=2\log (| \cos \pi(x - b-1/2)|)=2\log (|\sin \pi (x-b)|)\\
 \overline{\psi}^b(x)&:=\psi^b(x)-\mu(\psi^b)
\end{align*}
instead of the function 
$\phi^b=2\log (|x-b|)$ we consider. 
On the one hand this will simplify the calculations as we can use some statements from \cite{GKS}. On the other hand, the potential $\psi^b$ is of interest by itself, see Remark \ref{rem: Thue Morse pot}.
Later in Lemma \ref{lem:cylinder-to-orbit} we will show that the statements can easily be transferred.

As our observable has a singularity only to one side, the one-sided large deviation can easily be proven and is summarized in the following lemma:  
\begin{lemma}[One sided large deviation]\label{lem:upper_bound}
There are functions $g,C:\mathbb{R}_{>0}\to \mathbb{R}_{>0}$ (independent of $n$) such
that 
$$\mu(S_n(\overline{\psi}^b) \geqslant n\varepsilon)\leqslant C(\varepsilon)  e^{-ng(\varepsilon)},$$ 
for all $n\geqslant 0$. 
\end{lemma}

In order to prove the other direction we use a method adopted from \cite[Appendix B]{BDT}.
For doing so, we first define the topological pressure for an observable with a singularity on the negative side, introduced as in \cite{GKS}.
For a number of calculations and presentations it will be useful to identify $x\in \TT$ with its binary representation $x = \sum_{i=1}^{\infty} x_i 2^{-i}$, where each $x_i \in \Sigma = \{0,1\}$ (being unique for all $x$ except for the countable set of dyadic points) and the doubling map corresponds to the shift action, which we again denote by $T$,
\[
T \colon \Sigma^{\N} \to \Sigma^{\N}, \quad
(x_i)_{i \in \N} \mapsto (x_{i+1})_{i \in \N}.
\]
We refer to the elements of $\Sigma^n$ as \emph{words of length $n$} (for $n \in \N_{0}$) and define for $\omega = \omega_1 \cdots \omega_n \in \Sigma^n$ the corresponding \emph{cylinder set} by
\[
\langle \omega \rangle \coloneqq \left\{(x_i)_{i \in \N} \in \Sigma^\N : x_1 \cdots x_n = \omega_1 \cdots \omega_n \right\},
\]
in particular, if $\omega\in \Sigma^0$, we have the empty word and $\langle \omega \rangle=\Sigma^\N$. Moreover, sometimes we write $w_{[u,v]}$, $v\geqslant u$, by which we mean the word $w_uw_{u+1}\ldots w_v$.

We define the topological pressure (when the limit exists)  for a measurable function $f:[0,1]\to\mathbb{R}$ that is bounded from above as 
\begin{align}
  \mathcal{P}_{\operatorname{top}} (t f) \, 
  &= \,\lim_{n\to\infty}\frac{1}{n}\log \left(\sum_{\omega\in \Sigma^n}
  \exp \bigl(t \sup_{x\in \langle \omega\rangle}S_n(f)(x)\bigr)\right)\notag\\
  &= \,\begin{cases}\lim_{n\to\infty} \frac{1}{n}\log\left(\sum_{\omega\in \Sigma^n}
  \sup_{x\in \langle \omega\rangle}\exp \bigl(t S_n(f)(x)\bigr)\right)& t\geqslant 0\\
  \lim_{n\to\infty} \frac{1}{n}\log\left(\sum_{\omega\in \Sigma^n}
  \inf_{x\in \langle \omega\rangle}\exp \bigl(t S_n(f)(x)\bigr)\right)&t<0.
  \end{cases}\label{top_pressure}
\end{align}
Note that $\mathcal{P}_{\operatorname{top}} (0)=\log (2)$, which is the metric entropy for the doubling map.

Then, the next main step in the proof of Theorem \ref{LDP} is given in the following proposition:
\begin{prop}[Differentiability of the pressure function]\label{prop: pressure diff}
 Write 
 $$p(t)=\mathcal{P}_{\operatorname{top}}(t(\psi^b+2\log(2)))=\mathcal{P}_{\operatorname{top}}(t\overline{\psi}^b).$$ If $b$ is preperiodic but not periodic with respect to $T$ we have 
\begin{align}
 p '(0-)&=\lim_{t\downarrow 0} \frac{p(-t)-p(0)}{t}\notag\\
 &= \lim_{t\downarrow 0} \lim_{n\to\infty} \frac{\log \left(\sum_{\omega\in \Sigma^n} \inf_{y\in \langle\omega\rangle}\exp \left( -t(S_n(\psi^b)(y)+n2\log(2))\right)\right)}{n t}-\frac{\log(2)}{t}
=0.\label{eq: pressure diff}
\end{align}
\end{prop}

The next step in order to obtain a large deviation estimate is to get this estimate for an averaged sum which we introduce as follows:
For a measurable function $f:[0,1]\to\mathbb{R}$ that is bounded from above we define
\[
S_n^I (f) (x)=\sup_{y\in \langle x_{[1,n]}\rangle} S_n (f)(y),
\]
where $x_{[1,n]}$ denotes the word $x_1\ldots x_n$. 
Hence, $S_n^I (f)$ is constant on the cylinder sets $\langle \omega\rangle$ for all $\omega\in \Sigma^n$.

For $S_n^I$ we obtain the following lemma:
\begin{lemma}[Bounds for $S_n^I (x)$ from below via the pressure function]\label{inf:bound}
Under the condition that $p(t) =o(t)$, as $t\uparrow 0$, there exists an interval $U=(0,u)$ and functions $C_1,g_1:U\to \mathbb{R}_{>0}$ (independent of $n$) such that for all $n\geqslant 0$, we have 
\[
\mu \left(x:   S_n^I (\overline{\psi}^b) (x)  \leqslant -n\varepsilon \right)\leqslant C_1(\varepsilon) e^{-n g_1 (\varepsilon)}.
\]
\end{lemma}

In the next steps, we will show that the exponential large deviation bounds of Lemma~\ref{inf:bound} for 
$S_n^I (\overline{\psi}^b)$ implies  exponential large deviations for $S_n (\overline{\psi}^b)$. 
(The corresponding upper-tail bound
$\mu\bigl(S_n(\overline\psi^b)\geqslant n\varepsilon\bigr)$
is already given by Lemma~\ref{lem:upper_bound}.)
More precisely, we have the following lemma:
\begin{lemma}[Exponential bounds for $S_n(\overline{\psi}^b)$]\label{lem:cylinder-to-orbit}
Assume there exists an interval $U=(0,u)$ and functions $C_1,g_1:U\to\mathbb{R}_{>0}$
such that for all $n\geqslant 0$ we have
\[
\mu ( x :   S_n^I (\overline{\psi}^b)(x)  \leqslant -n \varepsilon )\leqslant C_ 1(\varepsilon) e^{-n g_1(\varepsilon)},
\] 
then there exist functions $C_2,g_2:U\to\mathbb{R}_{>0}$ (independent of $n$) such that for all $n\in\mathbb{N}$ we have
\[
\mu (x :   S_n (\overline{\psi}^b)  (x)   \leqslant -n \varepsilon )\leqslant C_2(\varepsilon)  e^{-n g_2 (\varepsilon)}.
\] 
\end{lemma}
 Combining these two estimates shows that
\(\overline\psi^b\) has exponential decay of deviations.

Finally, we compare $\phi^b$ with $\psi^b$ in the following lemma: 

\begin{lemma}[Comparing the large deviation of $\phi^b$ and $\psi^b$ ]\label{lem:psi-to-phi}
Let $\varepsilon >0$. Suppose there exists a neighborhood $U_1$ of $0$ such that for 
all  $\varepsilon \in U_1\backslash\{0\}$,  there is a constant $I_1(\varepsilon)>0$ and $C_1>0$ such
that $\mu(|S_n(\psi^b-\mu(\psi^b) )| \geqslant n\varepsilon)\leqslant C_1 e^{-nI_1(\varepsilon)} $ for all $n\geqslant 0$. Then there exists a $C_2>0$ and constant  $I_2(\varepsilon)>0$
defined on a neighborhood $U_2$ of $0$ such that for all $\varepsilon\in U_2\backslash\{0\}$ and all $n\geqslant 0$ we have
$$\mu(
|S_n(\phi^b-\mu (\phi^b) )| \geqslant n\varepsilon) <C_2 e^{-n I_2 (\varepsilon)}.$$
\end{lemma}

Combining the statements of the previous lemmas gives the statement of the theorem.

In the next section (Section \ref{sec: diff pressure}) we will prove Proposition \ref{prop: pressure diff}. 
After that, in Section \ref{sec: proofs large dev} we will prove all other statements of this section giving a complete proof of Theorem~\ref{LDP}. However, before proceeding with proving the just given statements we close this section with a remark regarding the rate function. 

\begin{remark}
Note that our methods don't allow for a rate function, i.e.\ we don't make a claim about the convexity of $I$ in \eqref{eq: main ineq}.
Indeed, examples exist in the literature~\cite{Bra89,OP88,BS93,Chu11} of
stationary processes which have exponential large deviations but a rate
function does not exist. 

In particular there is an example of a mean zero
bounded function $f$ taking only $3$ values on an aperiodic recurrent
Markov chain $(X_n)$ with a countable state space such that the system has
exponential large deviations but does not have a rate function. In this example defining $S_n=\sum_{j=0}^{n-1} f(X_j)$ for all $\varepsilon >0$, there
exist constants $C(\varepsilon)$, $0<\gamma<1$ such that
$\mathbb{P}( |\frac{S_n}{n}|>\varepsilon)\leqslant C(\varepsilon)e^{-\gamma n}$, giving
exponential convergence in the strong law of large numbers yet there is no
rate function controlling the rate of decay. However
in these examples $f$ is a coboundary~\cite{NT_LD} and there exists
$\psi \in L^2 (\mathbb{P})$ such that $f(X_j)=\psi (X_{j+1})-\psi (X_j)$
for all $j\geqslant 0$.  Hence  the example is degenerate (the variance
$\sigma^2=0$). 

We will see later (proof of Propostion~\ref{prop: pressure diff}) that $\psi^{1/2}$ is a coboundary as well (while $\psi^{b}$ for other values of $b$ is not). Hence, in the first case we also would not expect to obtain a rate function. 
\end{remark}

\section{Differentiability of the pressure function - Proof of Proposition \ref{prop: pressure diff}}\label{sec: diff pressure}
The statement of Proposition \ref{prop: pressure diff} is equivalent to the one-sided differentiability of the pressure function at $0$.
For a function $f$ for which the pressure function is defined we have
\begin{align}
 p_f'(0-)&\leqslant\lim_{t\downarrow 0} \frac{\mathcal{P}_{\operatorname{top}} (-t f)-\mathcal{P}_{\operatorname{top}} (0)}{t}\notag\\
 &= \lim_{t\downarrow 0} \lim_{n\to\infty} \frac{\log \left(\sum_{\omega\in \Sigma^n} \inf_{y\in \langle\omega\rangle}\exp \left( -t(S_n(f)(y))\right)\right)}{n t}-\frac{\log(2)}{t}. \label{eq: p' for Hoelder}
\end{align}
We will first describe how we would calculate the derivative of the pressure at $0$ for a H\"older continuous function $f$ with mean $0$ and show then that - given our potential is not H\"older continuous - the problem is indeed more complicated. By H\"older continuity of $f$ we immediately get 
\begin{align*}
\mathcal{P}_{\operatorname{top}} (t f)&= \lim_{n\to\infty} \frac{\log \left(\sum_{\omega\in \Sigma^n} \inf_{y\in \langle\omega\rangle}\exp \left( -t(S_n(f)(y))\right)\right)}{n}\\
&=\lim_{n\to\infty} \frac{\log \left(\sum_{\omega\in \Sigma^n} \sup_{y\in \langle\omega\rangle}\exp \left( -t(S_n(f)(y))\right)\right)}{n}.
\end{align*}
Given that
$n\mapsto a(2^{n})$ with 
\begin{align*}
 a(n)=\frac{\log \left(\sum_{\omega\in \Sigma^n} \sup_{y\in \langle\omega\rangle}\exp \left( -t(S_n(f)(y))\right)\right)}{n}
\end{align*}
is subadditive, using Fekete's lemma
it is possible to switch the limits in \eqref{eq: p' for Hoelder}. Using then the generalized mean  equality 
\[
\lim_{t\downarrow 0} \bigg(\frac{1}{N} \sum_{j=1}^N a_j^t\bigg)^{\frac{1}{t}}
=\bigg(\prod_{j=1}^N a_j\bigg)^{\frac{1}{N}}
\]
we get
\begin{align}
  p_f'(0-)
 &\leqslant  \lim_{n\to\infty}\lim_{t\downarrow 0} \frac{\log \left(\sum_{\omega\in \Sigma^n} \sup_{y\in \langle\omega\rangle}\exp \left( -t(S_n(f)(y))\right)\right)}{n t}-\frac{\log(2)}{t}\notag\\
  &=  \lim_{n\to\infty}\lim_{t\downarrow 0} \frac{\log \left(\sum_{\omega\in \Sigma^n} \sup_{y\in \langle\omega\rangle}\exp \left( -t(S_n(f)(y))- n\log(2)\right)\right)^{1/t}}{n }\notag\\
 &= \lim_{n\to\infty}
 \frac{ \log \left(\prod_{\omega\in \Sigma^n} \sup_{y\in \langle\omega\rangle} \exp\left( -(S_n(f)(y))\right)\right)^{2^{-n}}}{n}\notag\\
 &= \lim_{n\to\infty}
 \frac{ \sum_{\omega\in \Sigma^n} \sup_{y\in \langle\omega\rangle} \left( -(S_n(f)(y))\right)}{n2^n}\notag\\
 &\leqslant\lim_{n\to\infty}
 \frac{ \sum_{\omega\in \Sigma^n} \left(\inf_{y\in \langle\omega\rangle} \left( -(S_n(f)(y))\right)+W\right)}{n2^n}\notag\\ 
 &\leqslant \lim_{n\to\infty}\frac{ \int S_n(f) \,\mathrm{d}\lambda +W }{n}
 =\lim_{n\to\infty}\frac{W}{n}=0,
 \label{eq: limit switch Hoelder}
\end{align}
where $W$ denotes the bounded
distortion constant for $f$.

However, in our case having the unbounded potential $\psi^b$ instead of the H\"older continuous potential $f$ we can not make use of this method and given we need to use the infimum in the definition of the pressure function, we are in the superadditive regime. 

The problem to compute the pressure function at $0$ in a superadditive situation has recently attracted attention in a number of settings, e.g.\ \cite{rush} in the area of of matrix-cocycle potentials, under fairly general assumptions on the cocycle and generalizations in \cite{MVpreprint}. Indeed, there are also examples known where the pressure function for $t<0$ is continuous, but not differentiable, \cite{feng}. 

Already from this point of view, the differentiability of the pressure function of $\psi^b$ at $0$ is of independent interest. 
Moreover, the potential $\psi^b$ plays a special role in that it is  the logarithm of the unique $g$-function of the generalized Thue--Morse measure as outlined in Remark \ref{rem: Thue Morse pot}. 

In the sequel, we will give the proof of Proposition \ref{prop: pressure diff} in two steps: first for $b=1/2$ in which case $\psi^b$ is a co-boundary, the pressure function thus becomes linear and secondly for $b\neq 1/2$. 
In the latter case the main idea is to first restrict ourselves to a subshift of $\mathbb{X}$ on which $\psi^b$ is H\"older continuous, then switch the limits as in \eqref{eq: limit switch Hoelder} and, in a technical calculation, carefully estimate the error terms that occur by not being on the full shift, but a subshift.

\begin{remark}\label{rem: Thue Morse pot}
 The potential $\psi^b$ is of particular interest as it is the logarithm of the unique $g$-function of the generalized Thue--Morse measure given by the Riesz product measure $\mu_b \coloneqq\lim_{N\to\infty}
 \prod_{k=0}^{N-1}\big(1+\cos(2\pi(2^k x-b-1/2))\big)\,\lambda$, where the limit is understood as the weak limit of probability measures on the torus $\mathbb{T}$, the finite
products are their densities with respect to normalized Lebesgue measure, see \cite{GKS}. (Note, the notation is slightly different here, the $\psi^c$ in there refers to $\psi^b$ with $b=c+1/2\mod 1$ here.) This Riesz product measure is the diffraction measure of an aperiodic system generated by a substitution dynamical system generated by a finite (in case $b$ is rational) or infinite (in case $b$ is irrational) alphabet. 
With the help of the pressure function it was possible to determine further properties of the measure $\mu_b$, in particular the local dimension spectrum or the $L^q$-spectrum. 

The potential $\psi^0$ corresponds to the classical Thue--Morse case. It is clear that our method here can not work as the pressure for negative $t$ equals $\infty$, see \cite{BGKS}. This is in line with Proposition~\ref{prop: periodic part}.
\end{remark}

We start by giving the proof for the special case $b=1/2$. 

\begin{proof}[Proof of Proposition \ref{prop: pressure diff} for $b=1/2$] 
This case is special as $\psi^{1/2}$ is a coboundary 
\begin{equation}
\label{EQ:coboundary}
\psi^{1/2}(x) = 2 \log (|\sin(2 \pi x)|) - 2 \log (|\sin(\pi x)|) - 2 \log(2),
\end{equation}
for all $x \neq 0$ which results in the special representation of the pressure
$\mc P_{\operatorname{top}}(t \psi^{1/2})  = \max\{ 1-2t, 0 \}\log(2)$, see \cite{fan_generalized_2020}.
Hence, we obtain 
 \begin{align*}
   p '(0-)&\leqslant \lim_{t\downarrow 0} \frac{p(-t)-p(0)}{t}
 = \lim_{t\downarrow 0} \frac{(1+2t)\log (2)-2\log(2)t-\log (2)}{t}
 = 0,
\end{align*}
giving the statement of the proposition.
\end{proof}

The proof of $b\neq 1/2$ is substantially more involved. In the following we will give the main steps of the proof and postpone the proofs of the technical lemmas needed to the following Subsection \ref{subsec: proofs pressure}. 
\begin{proof}[Proof of Proposition \ref{prop: pressure diff} for $b\neq 1/2$ (main steps)]
In here, as stated in the proposition, we will always assume that $b$ is preperiodic but not periodic with respect to $T$. During the proof we will often distinguish between $b$ a dyadic point, i.e.\ the binary representation of $b$ ending on $0^{\infty}$ or $1^{\infty}$ and all other cases. However, the main idea of the proof remains the same.

We will first introduce some additional notation that is required for the proof and then give the main lemmas needed to show the proposition. 

If $b\neq 0$ is a dyadic point, then it can be written as $u10^{\infty}=u01^{\infty}$. As the case $b=1/2$ has already been studied we assume that $u$ is not the empty word. 
Moreover, for the following let 
$b=u10^{\infty}$ and 
$b^{\ast}=u01^{\infty}$ be the alternative representation of $b$.

If $b\neq 0$ is not a dyadic point, then it can be written as $uv^{\infty}$, where $|v|=:\mathsf{p}>1$ and the convention that $u$ is the shortest word such that $b$ can be written as $uv^{\infty}$. As we assume that $b$ is not periodic, we have that $|u|\geqslant 1$. 
Moreover, we set $\mathcal{V}=\{v_{[j,\mathsf{p}]}v_{[1,j-1]}, j\in \{1,\ldots,\mathsf{p}\}\}$. Here, $v_{[j,\mathsf{p}]}v_{[1,j-1]}$ denotes the concatenation of the two words $v_{[j,\mathsf{p}]}$ and $v_{[1,j-1]}$.
For the following we will always assume that
$m> 2|u|$. 
Sometimes we will use $b^{\ast}$ also in case $b$ is not a dyadic point in which case we formally identify $b=b^{\ast}$. 

In the next steps we will construct a sequence of subspaces $\widetilde{\mathbb{X}}^b_m$ that avoids the singularity at $b$ and we will compare the pressure on $\widetilde{\mathbb{X}}^b_m$ with the pressure on $\mathbb{X}$. 

For doing so, we first let $\widetilde{\Sigma}^{mn}_m$ be the set of $mn$-letter words where for each $w\in \widetilde{\Sigma}^{mn}_m$ and each $j\leqslant n$ we have 
\begin{align*}
 w_{[(j-1)m+1, jm]}
 &\notin\begin{cases}
         \{0^{m},1^m\}&\text{if }b\text{ is a dyadic point}\\
         \{v^{m/\mathsf{p}}, v\in\mathcal{V}\}&\text{otherwise.}
        \end{cases}
\end{align*}

One can think of $\widetilde{\Sigma}^{mn}_m$ as a word of $n$ letters of a new alphabet that consists of all $m$-letter words of the alphabet $\{0,1\}^m\backslash\{0^m,1^m\}$ in case of a dyadic point and  $\{0,1\}^m\backslash \{v^{m/\mathsf{p}},v\in\mathcal{V}\}$ otherwise.  We will denote the new shift space by $\widetilde{\X}_m^b$.

 Moreover, we define a function $\theta_m:\Sigma^{mn}\to\widetilde{\Sigma}_{m}^{mn}$ in the following way: If $b$ is a dyadic point, for $w\in \Sigma^m$ we set
\begin{align}
 \theta_m(w)=
 \begin{cases}
              w &\text{if }w\notin \{0^m, 1^m\}\\
              u10^{m-|u|-1}&\text{if }w=0^m\\
              u01^{m-|u|-1}&\text{if }w=1^m.
             \end{cases}\label{eq: thetam precon}
\end{align}
If $b$ is not a dyadic point, for $w\in \Sigma^m$ we set
\begin{align}
 \theta_m(w)=\begin{cases}
             w&\text{if }w\notin\{v^{m/\mathsf{p}}: v\in\mathcal{V}\}\\
             0^juv^{(m-j-|u|)/\mathsf{p}}&\parbox[t]{1\textwidth}{
             if $w=v_{[\ell,\mathsf{p}]}v^{(m-\mathsf{p}+\ell-1)/\mathsf{p}}$,\\  where $\ell\equiv |u|+j+1\mod \mathsf{p}$ and $j\in\{0,\ldots, \mathsf{p}-1\}$.}
             \end{cases}\label{eq: thetam preper}
\end{align}
We set 
$\theta_m(\Sigma^{im}\Sigma^{jm})=\theta_m(\Sigma^{im})\theta_m(\Sigma^{jm})$, where $\Sigma^{im}\Sigma^{jm}$ denotes the concatenation of two words of length $im$ and $jm$ respectively. 
 Note that in the dyadic case for $v\in \widetilde{\Sigma}^m_m$ we have that $\#\theta_m^{-1}(v)=2$ if and only if $v\in \{u01^{m-|u|-1}, u10^{m-|u|-1}\}$.

Finally, note that for any $w\in \Sigma^{nm}$ we have 
that at most $2m-1$ consecutive letters in $\theta_m(w)$ coincide with the first letters of $b$.

The main idea of the proof of Proposition \ref{prop: pressure diff} is to determine what is happening on the subshift $\widetilde{\mathbb{X}}^{b}_m$ on which we are in the classical setting of having a H\"older continuous potential and after taking care of a number of error terms we will be able use a similar approach as demonstrated in \eqref{eq: limit switch Hoelder}
Then we let $m$ tend to infinity. 
To start with these considerations we first need the following proposition comparing $S_n(\overline{\psi}^b)$ on the fullshift and on the subshift.

\begin{prop}[Comparing $S_n(\overline{\psi}^b)$ on $\langle w\rangle$ and $\langle \theta_m(w)  \rangle \cap \widetilde{\mathbb{X}}_{m}^b$]
\label{PROP:word-regularisation}
For all $\varepsilon > 0$, there exists $M\in \mathbb{N}$ such that for all $m\geqslant M$ there exists $N(m)$ such that for all $n\geqslant N(m)$ and all $w \in \Sigma^{mn}$ we have 
\[
\sup_{x \in \langle w \rangle} S_{nm}(\overline{\psi}^b)(x)
\geqslant
\sup_{y \in \langle \theta_m(w)  \rangle \cap \widetilde{\mathbb{X}}_{m}^b } S_{nm}(\overline{\psi}^b)(y) - \varepsilon n m.
\]
\end{prop}

 By Proposition \ref{PROP:word-regularisation} we have
\begin{align*}
 p '(0-)&\leqslant \lim_{t\downarrow 0} \frac{p(-t)-p(0)}{t}\\
 &= \lim_{t\downarrow 0} \lim_{n\to\infty} \frac{\log \left(\sum_{\omega\in \Sigma^n} \inf_{y\in \langle\omega\rangle}\exp \left( -tS_n(\overline{\psi}^b)(y)\right)\right)}{n t}-\frac{\log(2)}{t}\\
  &\leqslant  \lim_{t\downarrow 0} \lim_{n\to\infty} \frac{\log \left(\sum_{\omega\in \Sigma^{nm}} \inf_{y\in \langle\omega\rangle}\exp \left(-tS_{nm}(\overline{\psi}^b)(y)\right)\right)}{nm t}-\frac{\log(2)}{t}\\
  &\leqslant  \lim_{t\downarrow 0} \lim_{n\to\infty} \frac{\log \left(\sum_{\omega\in \Sigma^{nm}} \sup_{y\in \langle\theta_m(\omega)\rangle\cap\widetilde{\X}_m^b}\exp \left( -t(S_{nm}(\overline{\psi}^b)(y)-nm\varepsilon)\right)\right)}{nm t}-\frac{\log(2)}{t}.
\end{align*}

For fixed $t$, let 
\begin{align*}
 \widetilde{F}(n,t)= \frac{\log \left(\sum_{\omega\in \Sigma^{nm}} \sup_{y\in \langle\theta_m(\omega)\rangle\cap\widetilde{\X}_m^b}\exp \left( -t(S_{nm}(\overline{\psi}^b)(y)-nm\varepsilon)\right)\right)}{nm t}
\end{align*}
and set $F(n,t)= \widetilde{F}(m^n,t)$. 
Then, we have 
\begin{align*}
 m^{n+1} F(n+1,t)= m^{n+1} \widetilde{F}(m^{n+1},t)\leqslant m^{n} (m\widetilde{F}(m^n,t))= m^{n+1} F(n,t). 
\end{align*}
Hence, $F(\cdot, t)$ is a decreasing function. On the other hand, using the properties of the generalized mean, we have that $F(n, \cdot )$ is increasing. Hence, we can exchange the limits to obtain
\begin{align}
 p' (0-)&\leqslant \lim_{t\downarrow 0} \lim_{n\to \infty}
 \frac{\log \left(\sum_{\omega\in \Sigma^{nm}} \sup_{y\in \langle\theta_m(\omega)\rangle\cap\widetilde{\X}_m^b}\exp \left( -t(S_{nm}(\overline{\psi}^b)(y)-nm\varepsilon)\right)\right)}{nm t}-\frac{\log(2)}{t}\notag\\
 &=  \lim_{n\to \infty}\lim_{t\downarrow 0}
 \frac{\log \left(\sum_{\omega\in \Sigma^{nm}} \sup_{y\in \langle\theta_m(\omega)\rangle\cap\widetilde{\X}_m^b}\exp \left( -t(S_{nm}(\overline{\psi}^b)(y)-nm\varepsilon)\right)\right)}{nm t}-\frac{\log(2)}{t}\notag\\ 
  &=  \lim_{n\to \infty}\lim_{t\downarrow 0}
 \frac{1}{nm} \log \left( 2^{-nm}\sum_{\omega\in \Sigma^{nm}} \sup_{y\in \langle\theta_m(\omega)\rangle\cap\widetilde{\X}_m^b}\exp \left( -(S_{nm}(\overline{\psi}^b)(y)-nm\varepsilon)\right)^t\right)^{1/t}.\notag
 \end{align}
 Using then the properties of the generalized mean yields 
 \begin{align}
   p' (0-)&\leqslant  \lim_{n\to\infty}
 \frac{1}{nm} \log \left(\prod_{\omega\in \Sigma^{nm}} \sup_{y\in \langle\theta_m(\omega)\rangle\cap\widetilde{\X}_m^b}\exp \left( -S_{nm}(\overline{\psi}^b)(y)+nm\varepsilon\right)\right)^{ 2^{-nm}}\notag\\ 
    &=  \lim_{n\to\infty}
 \frac{1}{nm 2^{nm}} \sum_{\omega\in \Sigma^{nm}} \sup_{y\in \langle\theta_m(\omega)\rangle\cap\widetilde{\X}_m^b} -S_{nm}(\overline{\psi}^b)(y)+nm\varepsilon\notag\\ 
      &=  \varepsilon-\lim_{n\to\infty}
  \frac{1}{nm 2^{nm}} 
  \sum_{\omega\in \Sigma^{nm}} \inf_{y\in \langle\theta_m(\omega)\rangle\cap\widetilde{\X}_m^b}
  S_{nm}(\overline{\psi}^b)(y),\label{eq: p(0-)}
\end{align}  
for all $\varepsilon>0$.

Our next goal is to compare $\sum_{w\in \Sigma^{nm}}\inf_{y\in \langle\theta_m(\omega)\rangle\cap\widetilde{\X}_m^b}
  S_{nm}(\overline{\psi}^b)(y)$ with $\int
  S_{nm}(\overline{\psi}^b)\,\mathrm{d}\mu=0$.
However, an intermediate step will first be to show the following lemma. 
\begin{lemma}\label{lem: split up of sums}
 For all $\varepsilon>0$ there exists $M_0\in\mathbb{N}$ such that for all $m\geqslant M_0$ there exists $N(m)$ such that for all  $n\geqslant N(m)$ we have
\begin{align*}
\sum_{\omega\in \Sigma^{nm}} \sum_{i=1}^{nm} 
    \sup_{x\in \langle\theta_m(\omega)\rangle}\psi^b(T^{i-1}(x))- \inf_{y\in \langle\theta_m(\omega)\rangle\cap\widetilde{\X}_m^b}\psi^b(T^{i-1}(y))
    \leqslant \varepsilon nm2^{nm}.
\end{align*}
\end{lemma}
From Lemma \ref{lem: split up of sums} we then obtain 
\begin{align} 
 \MoveEqLeft-\lim_{n\to\infty}
  \frac{1}{nm 2^{nm}} 
  \sum_{\omega\in \Sigma^{nm}} \inf_{y\in \langle\theta_m(\omega)\rangle\cap\widetilde{\X}_m^b}
  S_{nm}(\overline{\psi}^b)(y)\notag\\
  &\leqslant -\lim_{n\to\infty}
  \frac{1}{nm} 
  \sum_{\omega\in \Sigma^{nm}} \int_{\langle\theta_m(\omega)\rangle}
  S_{nm}(\overline{\psi}^b)\,\mathrm{d}\mu+2\varepsilon\notag\\
  &\leqslant -\lim_{n\to\infty}
  \frac{1}{nm} 
  \sum_{\omega\in \Sigma^{nm}} \int_{\langle\omega\rangle}
  S_{nm}(\overline{\psi}^b)\,\mathrm{d}\mu+2\varepsilon+E(n,m)
  \leqslant 2\varepsilon+E(n,m),\label{eq: p(0-) 1 new}
\end{align}
where
\begin{align*}
 E(n,m)
 &=\frac{1}{nm}\sum_{\omega\in \Sigma^{nm}} \left(\int_{\langle\omega\rangle}S_{nm}(\overline{\psi}^b)\,\mathrm{d}\mu
 -
 \int_{\langle\theta_m(\omega)\rangle}S_{nm}(\overline{\psi}^b)\,\mathrm{d}\mu\right).
\end{align*} 
The final step is to estimate the error term $E(n,m)$; the proof of the following lemma we postpone to the following subsection (Subsection \ref{subsec: proofs pressure}).
\begin{lemma}\label{lem: Enm estim}
 For all $\varepsilon>0$ there exists $M\in\mathbb{N}$ such that for all $m\geq M$ and all $n$ we have 
 \begin{align*}
  E(n,m)\leq \varepsilon. 
\end{align*}
\end{lemma}

Combining then \eqref{eq: p(0-)} with \eqref{eq: p(0-) 1 new} and the just given error estimate yields 
$p'(0-)
\leqslant 4\varepsilon$, 
for $m$ sufficiently large. 
Since $\varepsilon$ was arbitrary, the statement of the proposition follows. 
\end{proof}

In the following subsection we will give the proofs of the just stated lemmas:
\subsection{Proofs concerning the differentiability of the pressure function}\label{subsec: proofs pressure}
In order to prove Proposition \ref{PROP:word-regularisation} we first introduce some additional notation and state some technical lemmas.

Next, we start by introducing some notation and citing some lemmas.
Again, we denote by $b_{[m,n]}$ the finite word of the $m$th digit to the $n$th digit of $b$.
In case $b_{[1,n]}= b^{\ast}_{[1,n]}$ we introduce
\begin{equation}
 \mc G_n = \mc G_n(b) = b_{[1,n]} + \{ 0,1,-1 \}; \quad \mc G =\mc G (b)= \bigcup_{n \geqslant 1} \mc G_n,\label{eq: Gn non dyad}
\end{equation}
where the $+$ is understood as addition of $b_{[1,n]}$ and $0,1,-1$ as numbers in the binary system 
together with the conventions $0^n - 1 = 1^n$ and $1^n + 1 = 0^n$.
That is, $\{\langle w \rangle\}_{w \in \mc G_n}$ consists of the length-$n$ cylinder set containing $b$, along with its right and left neighbour.
Moreover, we set 
\begin{equation}
 \mc G_n = \left\{b_{[1,n]}, b^{\ast}_{[1,n]} \right\},\label{eq: Gn dyad}
\end{equation}
for all $n$ with $b_n \neq b^{\ast}_n$,
without altering any of the decisive properties above. Moreover, we also set $\mc G_0$ as the empty word.
If $b$ is not a dyadic word, we will always use \eqref{eq: Gn non dyad} while in case that $b$ is dyadic, we formally use \eqref{eq: Gn non dyad} for small $n$ and \eqref{eq: Gn dyad} for $n\geqslant |u|+1$. However, as we are mostly interested in large $n$, we can morally think of \eqref{eq: Gn non dyad} as the definition of $\mc G_n$ non dyadic points $b$ and \eqref{eq: Gn dyad} for dyadic points $b$.

\begin{lemma}[{\cite[Lemma 9.11]{GKS}}]
\label{LEM:psi-distance-bound}
For every $x$, we have
\[
2 \log (2 \rho(x,b)) \leqslant \psi^b(x)
\leqslant 2 \log (\pi \rho(x, b)).
\]
Furthermore, if $x,y \in \mathbb{T} \setminus B_r(b)$ for some $r>0$, we have the Lipschitz estimate
\[
|\psi^b(x) - \psi^b(y)| \leqslant \frac{2}{r} \rho(x,y).
\]
\end{lemma}

The next lemma is a reformulation of \cite[Lemma 9.12]{GKS} but can also easily be deduced from Lemma \ref{LEM:psi-distance-bound}.
\begin{lemma}
\label{LEM:psi-on-hitting}
Let $w \in \mc G_n$ and $x\in\langle w\rangle$.
Then we have 
\[
\psi^b(x) \leqslant 
(-n +1)\, 2 \log(2).
\]
On the other hand, if $w\notin \mc G_{n+1}$ and $x\in\langle w\rangle$, then we have
\[
\psi^b(x)\geqslant -n 2 \log(2).
\]

If $w_{[1,n]}\in \mc G_n$ but $w_{[1,n+1]}\notin \mc G_{n+1}$ and $x,y\in\langle w\rangle$, then we have
\[
|\psi^b(x) - \psi^b(y)|
\leqslant  2^{-n + 3}.
\]
\end{lemma}

Given a word $w$ and $j \in \N$ with $j < |w|$, the \emph{$j$-hitting times} (with respect to $w$) are given by
\[
\mc H_j(w) = \left\{ k \leqslant j: w_{[k,j]} \in \mc G \mbox{ and } w_{[k,j+1]} \notin \mc G \right\}.
\]
Likewise, the \emph{(full) hitting times} (with respect to $w$) are
\[
\mc H(w) = \left\{ k \leqslant |w|: w_{[k,|w|]} \in \mc G \right\}.
\]

\begin{proof}[Proof of Proposition~\ref{PROP:word-regularisation}]
Before starting with the technical details of the proof, we will first give an idea of the construction. Let $\widetilde{x}$ be the point such that $\sup_{x\in \langle w\rangle}S_n(\psi^b)(x)=S_n(\psi^b)(\widetilde{x})$.
Loosely speaking, if $T^j(\tilde{x})$ is close to the singularity for some $j\in \mathbb{N}$, then necessarily $w_{[j+1,j+r]}\in \{b_{[1,r]}, b^{\ast}_{[1,r]}\}$ for some large $r$. The idea of the construction is to not have one single very negative summand $\psi^b(T^j(x))$ but a number of negative summands $\psi^b(T^{q+im}(\theta_m(\tilde{x}))$ where $q=\min\{v\in\mathbb{N}: v\geqslant j\text{ and }m\mid v+1\}$ and $i\in \{0,\ldots, \lfloor r/m\rfloor\}$. All of those summands are smaller in modulus than $\psi^b(T^j(x))$ but add up to a comparable size as $\psi^b(T^j(x))$.

Next, we introduce some additional notation we need during the proof. Let $w\in \Sigma^{mn}$ with $m\geqslant |u|+2$ and let $w'=\theta_m(w)$. Moreover, we set $x'=x'(w)\in \langle w \rangle$ as follows: 
\begin{align}\label{eq: def x'}
 x'=\begin{cases}
     w(01)^{\infty}&\text{if }b\text{ is a dyadic point}\\
     w0^{\infty}&\text{otherwise.}
    \end{cases}
\end{align}

Note that we immediately have 
$\sup_{x \in \langle w \rangle} S_n(\psi^b)(x)\geqslant S_n(\psi^b)(x')$.
In the following we will compare 
 $S_n(\psi^b)(x')$ with $\sup_{y \in \langle w'  \rangle \cap \widetilde{\mathbb{X}}_{m}^b } S_n(\psi^b)(y)$ where we denote the point on which the supremum is attained by $y'$.

We will separately estimate $\sum_{k=jm+1}^{(j+1)m} \psi^b (T^{k-1}(x'))$, for different values of $j$. Here we distinguish different cases. First, for the following cases, we will assume $j\leqslant n-1$.

In the following, to improve readability,  for each case we first give a proof for $b$ being a dyadic point and then give a shortened proof for the remaining case.

\noindent\textbf{1)} $w_{[jm+1, (j+1)m]}=w'_{[jm+1, (j+1)m]}$ and 
$\mc H(w_{[jm+1, (j+2)m]})\cap \N_{\leqslant m}=\emptyset$

\begin{adjustwidth}{.8cm}{0cm}
 \noindent\emph{Case 1: $b$ is a dyadic point:}
 
 This case is either obtained if $w_{[jm+1, (j+2)m]}=w'_{[jm+1, (j+2)m]}$ or 
 if $w_{[jm+1, (j+1)m]}=w'_{[jm+1, (j+1)m]}$ and $w_{[(j+1)m+1, (j+2)m]}\neq w'_{[(j+1)m+1, (j+2)m]}$ but the concatenation of the two words $w_{[jm+1, (j+2)m]}$ and $w'_{[jm+1, (j+2)m]}$ does not end on $b_{[1,q]}$ or $b^{\ast}_{[1,q]}$ with $q>m$.

 Let $A_{j,m}=\{jm+1,\ldots, (j+1)m\}$ 
 and set $\overline{A}_{j,m}=\{k \in A_{j,m}: x'_{[k,k+|u|+1]}\in \{u01, u10\}\}$ and $\underline{A}_{j,m}=A_{j,m}\backslash \overline{A}_{j,m}$. 
 First note that for any $k\in \underline{A}_{j,m}$ we have by Lemma \ref{LEM:psi-distance-bound} that 
 \begin{equation}
 |\psi^b(T^{k-1}(x'))-\psi^b(T^{k-1}(y'))|
 \leqslant \frac{2}{2^{-|u|-2}} 2^{k-(j+1)m}.\label{eq: diff x thetax equal0}
 \end{equation}

 Hence, 
 \begin{align}
  \sum_{k\in \underline{A}_{j,m}}|\psi^b(T^{k-1}(x'))-\psi^b(T^{k-1}(y'))|
  &\leqslant \sum_{k\in \underline{A}_{j,m}}2^{k+|u|+3-(j+1)m}\notag\\
  &\leqslant \sum_{k\in A_{j,m}}2^{k+|u|+3-(j+1)m}
  \leqslant 2^{|u|+4}.\label{eq: diff x thetax equal1}
 \end{align}
 Next, we analyse $k\in \overline{A}_{j,m}$ and let the elements of $\overline{A}_{j,m}$ be denoted by $k_1< \ldots <k_i$ with $i=\#\overline{A}_{j,m}$. Then we have for all $\ell<i$ that $T^{k_\ell-1}(x'), T^{k_\ell-1}(y')\in \T\backslash B_r(b)$ with $r\geqslant 2^{k_{\ell}-k_{\ell+1}-|u|-2}$. This can be seen as follows: Latest at the positions $k_{\ell+1}+|u|, k_{\ell+1}+|u|+1$ we have two consecutive letters that are different. 
 Hence, by Lemma \ref{LEM:psi-distance-bound} we obtain 
 \begin{align}
  |\psi^b(T^{k_\ell-1}(x'))-\psi^b(T^{k_\ell-1}(y'))|
  \leqslant \frac{2}{2^{k_{\ell}-k_{\ell+1}+|u|-2}} 2^{-(j+1)m+ k_\ell}
  =2^{-(j+1)m+|u|+3+k_{\ell+1}}. \label{eq: diff x thetax equal1a}
 \end{align}
 Summing over $\ell$ we obtain 
 \begin{align}
  \sum_{\ell=1}^{i-1}|\psi^b(T^{k_\ell-1}(x'))-\psi^b(T^{k_\ell-1}(y'))|
  &\leqslant \sum_{\ell=1}^{i-1}2^{-(j+1)m+|u|+3+k_{\ell+1}}\notag\\
  &\leqslant \sum_{\ell'=jm+2}^{(j+1)m}2^{-(j+1)m+|u|+3+\ell'}
  \leqslant 2^{|u|+4}.\label{eq: diff x thetax equal2}
 \end{align}

 Next, we consider $k_i$. Here, we distinguish the two subcases. 
 First, let us assume that $w_{[jm+1, (j+2)m]}=w'_{[jm+1, (j+2)m]}$ holds. This implies $\rho(T^{k_i-1}(x'),T^{k_i-1}(y'))\leqslant 2^{-(j+2)m+ k_i-1}$.
 By the above assumption $\mc H(w_{[jm+1, (j+2)m]})\cap \N_{\leqslant m}=\emptyset$ we have that $T^{k_i-1}(x'),T^{k_i-1}(y')\in \T\backslash B_r(b)$ with $r\geqslant 2^{k_i-(j+2)m-1}$. 
 Thus, by Lemma \ref{LEM:psi-distance-bound} we obtain 
 \begin{align}
  |\psi^b(T^{k_i-1}(x'))-\psi^b(T^{k_i-1}(y'))|
  \leqslant \frac{2}{2^{k_i-(j+2)m-1}} 2^{-(j+2)m+ k_i-1}
  \leqslant 2.\label{eq: diff x thetax equal3}
 \end{align} 

 Next, let us assume 
 $$w_{[jm+1, (j+1)m]}=w'_{[jm+1, (j+1)m]}\quad \text{ and }\quad w_{[(j+1)m+1, (j+2)m]}\neq w'_{[(j+1)m+1, (j+2)m]}$$ but neither of the two words $w_{[jm+1, (j+2)m]}$ or $w'_{[jm+1, (j+2)m]}$ end on $b_{[1,q]}$ or $b^{\ast}_{[1,q]}$ with $q>m$.
 This in particular implies that $T^{k_i-1}(x'),T^{k_i-1}(y')\in \T\backslash B_r(b)$ with $r\geqslant 2^{k_i-(j+1)m+1}$.
 Thus, by Lemma \ref{LEM:psi-distance-bound} we obtain 
 \begin{align}
  |\psi^b(T^{k_i-1}(x'))-\psi^b(T^{k_i-1}(y'))|
  \leqslant \frac{2}{2^{k_i-(j+1)m+1}} 2^{-(j+1)m+ k_i}
  \leqslant 4.\label{eq: diff x thetax equal3a}
 \end{align}
Finally, combining \eqref{eq: diff x thetax equal1}, \eqref{eq: diff x thetax equal2}, \eqref{eq: diff x thetax equal3} and \eqref{eq: diff x thetax equal3a}
we obtain  
\begin{align}
   \sum_{k\in A_{j,m}}|\psi^b(T^{k-1}(x'))-\psi^b(T^{k-1}(y'))|
  &\leqslant 4+2^{|u|+5}=:U_1\label{eq: U1}
\end{align}
which is independent of $m$, $j$ and $n$ but depends on $b$. 
\newline
 \noindent\emph{Case 2: $b$ is not a dyadic point}
 
 The distinguished cases remain the same, however, we need to define $\overline{A}_{j,m}=\{k\in A_{j,m}: x'_{[k,k+|u|+1]}\in \{u\}$. \eqref{eq: diff x thetax equal0} and \eqref{eq: diff x thetax equal1} then change to 
 \begin{align*}
   |\psi^b(T^{k-1}(x'))-\psi^b(T^{k-1}(y'))|
 &\leqslant \frac{2}{2^{-|u|}} 2^{k-(j+1)m}\quad\text{ and}\\
  \sum_{k\in \underline{A}_{j,m}}|\psi^b(T^{k-1}(x'))-\psi^b(T^{k-1}(y'))|
  &\leqslant 2^{|u|+2}.
 \end{align*}
 
 With the same notation of $k_1< \ldots <k_i$ and $i=\#\overline{A}_{j,m}$ we then have for all $\ell<i$ that $T^{k_\ell-1}(x'), T^{k_\ell-1}(y')\in \T\backslash B_r(b)$ with $r\geqslant 2^{k_{\ell}-k_{\ell+1}-|u|-\mathsf{p}}$ which follows from the fact that the representation $b=uv^{\infty}$ is chosen with the shortest possible word $u$. Thus,  $b_{[k_{\ell},k_{\ell+1}+|u|+\mathsf{p}-1]}$ can not be a prefix of $b$ anymore. 
 
Then, \eqref{eq: diff x thetax equal1a} and \eqref{eq: diff x thetax equal2} change to
 \begin{align*}
  |\psi^b(T^{k_\ell-1}(x'))-\psi^b(T^{k_\ell-1}(y'))|
  &\leqslant 2^{-(j+1)m+|u|+\mathsf{p}+1+k_{\ell+1}}\\ 
  \sum_{\ell=1}^{i-1}|\psi^b(T^{k_\ell-1}(x'))-\psi^b(T^{k_\ell-1}(y'))|
  &\leqslant 2^{|u|+\mathsf{p}+2}=:U_1.\label{eq: diff x thetax equal2}
 \end{align*}

\end{adjustwidth}

\noindent\textbf{2)} $w_{[jm+1, (j+1)m]}\in \{0^m, 1^m\}$ (in case $b$ is a dyadic point, otherwise $w_{[jm+1, (j+1)m]}=v^{m/\mathsf{p}}$ with $v\in\mathcal{V}$) and $w_{[j(m+1)+1, (j+2)m]}$ arbitrary

\begin{adjustwidth}{.8cm}{0cm}
 \noindent\emph{Case 1: $b$ is a dyadic point:} 
 
We only treat the case $w_{[jm+1, (j+1)m]}=0^m$ as the other case can be treated analogously. 
We have that $\theta_m(w_{[jm+1, (j+1)m]})=u10^{m-|u|-1}$.

 For all $\iota\in \{0, \ldots, |u|+1\}$ there exists $\ell\leqslant |u|+2$ (including possibly $\ell=0$) such that 
 $w_{[jm+1+\iota, jm+1+\iota+\ell]}\notin \mc G_{\ell+1}$. 
 Lemma \ref{LEM:psi-on-hitting}
 then implies 
 \begin{equation}
   \psi^b(T^{jm+\iota} (x'))\geqslant (-|u|-2) 2 \log (2)\quad\text{ for all }\iota \in \{0, \ldots, |u|+1\}.\label{eq: wjm 0 a}
 \end{equation}
 
 On the other hand, we have $w'_{[jm+1, (j+1)m]}\in \mc G_{m}$.
 By Lemma \ref{LEM:psi-on-hitting} this implies 
\begin{equation}
       \psi^b(T^{jm} (y'))\leqslant  2(-m+1) \log(2). \label{eq: wjm 0 b}                                                                                                         
\end{equation}

For all $\iota\in \{1,\ldots, |u|\}$ we simply use the trivial bound  
\begin{equation}
 \psi^b(T^{jm+\iota} (y'))\leqslant 0\quad\text{ for all }\iota \in \{1, \ldots |u|\}.\label{eq: wjm 0 c}
\end{equation}
Finally, we look at $T^{jm+\iota}(x')$ and $T^{jm+\iota}(y')$ for $\iota>|u|$ and compare the terms directly. 
For all $\iota\in\{|u|+1,\ldots, m-1\}$ we have $\rho(T^{jm+\iota}(x'), T^{jm+\iota}(y'))\leqslant 2^{-m+\iota}$ and moreover, $\min\{\rho(T^{jm+\iota} (x'), b), \rho(T^{jm+\iota}(y'), b)\}\geqslant 2^{-|u|-1}$. 
Hence, Lemma \ref{LEM:psi-distance-bound} implies 
\begin{align}
 \sum_{\iota=|u|+1}^{m-1} |\psi^b(T^{jm+\iota}(x'))-\psi^b(T^{jm+\iota}(y'))|
 \leqslant \sum_{\iota=|u|+1}^{m-1}\frac{2}{2^{-|u|-1}} 2^{-m+\iota}
 \leqslant 2^{|u|+2}.\label{eq: wjm 0 d}
\end{align}
Combining \eqref{eq: wjm 0 a}, \eqref{eq: wjm 0 b}, \eqref{eq: wjm 0 c}, and \eqref{eq: wjm 0 d} we obtain 
\begin{align*}
 \MoveEqLeft\sum_{\iota=0}^{m-1} \left(\psi^b(T^{jm+\iota}(y'))-\psi^b(T^{jm+\iota}(x'))\right)\\
 &\leqslant (|u|+3)(|u|+2)2\log(2)+2(-m+1) \log(2)+2^{|u|+2}
 \leqslant U_2-2\log(2)m, 
\end{align*}
where $U_2$ is independent of $m$, $j$ and $n$ but depends on $b$. 
\newline 

 \noindent\emph{Case 2: $b$ is not a dyadic point}

We have that $\theta_m(w_{[jm+1, (j+1)m]})=0^i uv^{m/\mathsf{p}-i-|u|}$ with $i\in\{0,\ldots,\mathsf{p}-1\}$ depending on which $v\in\mathcal{V}$ we are looking at.

 For all $\iota\in \{0, \ldots, |u|+1\}$ there exists $\ell\leqslant |u|+2\mathsf{p}+1$ (including possibly $0$) such that 
 $w_{[jm+1+\iota, jm+1+\iota+\ell]}\notin \mc G_{\ell+1}$. 
\eqref{eq: wjm 0 a} then becomes 
 \begin{equation*}
   \psi^b(T^{jm+\iota} (x'))\geqslant (-|u|-2\mathsf{p}-1) 2 \log (2)\quad\text{ for all }\iota \in \{0, \ldots, |u|+1\}.
 \end{equation*}
 
 On the other hand, we have $w'_{[jm+i+1, (j+1)m]}\in \mc G_{m-i}$.
\eqref{eq: wjm 0 b} then becomes
\begin{equation*}
       \psi^b(T^{jm} (y'))\leqslant  2(-m+i+1) \log(2)\leqslant -2\log(2)m+2\log(2)\mathsf{p}.                                                                                                       
\end{equation*}

For all $\iota\in \{0,\ldots, |u|+\mathsf{p}\}\backslash \{i\}$ we then simply use the trivial bound  
\begin{equation*}
 \psi^b(T^{jm+\iota} (y'))\leqslant 0
\end{equation*}
giving the analog of \eqref{eq: wjm 0 c}. 

For $\iota>|u|+\mathsf{p}$ the calculations are analogous with $\rho(T^{jm+\iota}(x'), T^{jm+\iota}(y'))\leqslant 2^{-m+\iota}$ and $\min\{\rho(T^{jm+\iota} (x'), b), \rho(T^{jm+\iota}(y'), b)\}\geqslant 2^{-|u|-\mathsf{p}-1}$ and as an analog of \eqref{eq: wjm 0 d} we obtain 
\begin{align*}
 \sum_{\iota=|u|+1}^{m-1} |\psi^b(T^{jm+\iota}(x'))-\psi^b(T^{jm+\iota}(y'))|
 \leqslant \sum_{\iota=|u|+\mathsf{p}+1}^{m-1}\frac{2}{2^{-|u|-\mathsf{p}-1}} 2^{-m+\iota}
 \leqslant 2^{|u|+\mathsf{p}+2}
\end{align*}
and combining the above inequalities, again we obtain 
\begin{align*}
 \sum_{\iota=0}^{m-1} \psi^b(T^{jm+\iota}(y'))-\psi^b(T^{jm+\iota}(x'))
 &\leqslant  U_2-2\log(2)m 
\end{align*}
with $U_2$ independent of $m$, $j$ and $n$, but depends on $b$. 
  
\end{adjustwidth}

\noindent\textbf{3)} $w_{[jm+1, (j+1)m]}=w'_{[jm+1, (j+1)m]}$ and 
$\mc H(w_{[jm+1, (j+q+1)m]})\cap \N_{\leqslant m}\neq\emptyset$ with $q\in \{1,\ldots, n-j-1\}$ maximal.

\begin{adjustwidth}{.8cm}{0cm}
 \noindent\emph{Case 1: $b$ is a dyadic point:} 
 
 First note that $q\leqslant n-j-1$ being maximal implies that either $\mc H(w_{[jm+1, (j+q+2)m]})\cap \N_{\leqslant m}=\emptyset$ or that $(j+q+1)m=nm$.
 Moreover, note that we automatically have $w_{[j(m+1)+1, (j+q+1)m]}\in \{0^{qm}, 1^{qm}\}$. As it is clear that $\mc H(w_{[jm+1, (j+q+1)m]})$ can only have one element, we may assume that $\{k\}=\mc H(w_{[jm+1, (j+q+1)m]})$. First we estimate 
 $|\psi^b(T^{jm+\iota}(x'))-\psi^b(T^{jm+\iota}(y'))|$ for $\iota\in\{0,\ldots, m-1\}\backslash\{k\}$. 
 In this case we have 
 \begin{align*}
  \min\left\{\rho(T^{jm+\iota}(x'),b), \rho(T^{jm+\iota}(y'),b)\right\}&\geqslant 2^{-|u|-1}\quad\text{ and }\\
  \rho(T^{jm+\iota}(x'),T^{jm+\iota}(y'))&\leqslant 2^{-m+\iota}.
 \end{align*}
 Hence, Lemma \ref{LEM:psi-distance-bound} implies 
 \begin{align}
  \sum_{\iota\in\{0,\ldots, m-1\}\backslash\{k\}}|\psi^b(T^{jm+\iota}(x'))-\psi^b(T^{jm+\iota}(y'))|
  \leqslant \sum_{\iota\in\{0,\ldots, m-1\}\backslash\{k\}}\frac{2}{ 2^{-|u|-1}} 2^{-m+\iota}
  \leqslant 2^{|u|+3}. \label{eq: wjm long a}
 \end{align}
 For $\iota=k$, we have $w_{[jm+k,(j+q+1)m]}\notin\mc G_{(q+2)m-k+1}$
 and Lemma \ref{LEM:psi-on-hitting} implies that 
\begin{align}
 \psi^b(T^{jm+k}(x'))
 &\geqslant 2 \log (2)(k- (q+2)m-1).\label{eq: wjm long b}
\end{align}
On the other hand,  $y'_{[k+jm,(j+1)m]}\in \mc G_{m-k+1}$. Hence, by Lemma \ref{LEM:psi-on-hitting} we obtain 
\begin{align}
 \psi^b(T^{jm+k}(y'))\leqslant 2\log(2) (k-m+2).\label{eq: wjm long c}
\end{align}
Combining \eqref{eq: wjm long a}, \eqref{eq: wjm long b}, and \eqref{eq: wjm long c}, we obtain 
\begin{align}
 \MoveEqLeft\sum_{\iota=0}^{m-1}\psi^b(T^{jm+\iota}(y'))-\psi^b(T^{jm+\iota}(x'))\notag\\
 &\leqslant 2^{|u|+3}+ 2\log(2) (k-m+2)- 2 \log (2)(k- (q+2)m-1)\notag\\
 &\leqslant 2^{|u|+3}+6\log(2)+ 2\log(2) (q+1)m\notag\\
 &=:U_3+ 2\log(2)(q+1)m.\label{eq: wjm long final} 
\end{align}
\newline 

 \noindent\emph{Case 2: $b$ is not a dyadic point:}

 Instead of $w_{[j(m+1)+1, (j+q+1)m]}\in \{0^{qm}, 1^{qm}\}$ we have $w_{[j(m+1)+1, (j+q+1)m]}=v^{qm/\mathsf{p}}$, for some $v\in\mathcal{V}$.

 For $\iota\in\{0,\ldots, m-1\}\backslash\{k\}$ with $\{k\}=\mc H(w_{[jm+1, (j+q+1)m]})$ we have 
 \begin{align*}
  \min\left\{\rho(T^{jm+\iota}(x'),b), \rho(T^{jm+\iota}(y'),b)\right\}&\geqslant 2^{-|u|-\mathsf{p}-1}\quad\text{ and }\\
  \rho(T^{jm+\iota}(x'),T^{jm+\iota}(y'))&\leqslant 2^{-m+\iota}
 \end{align*}
 and as an analog of \eqref{eq: wjm long a} we obtain 
 \begin{align*}
  \sum_{\iota\in\{0,\ldots, m-1\}\backslash\{k\}}|\psi^b(T^{jm+\iota}(x'))-\psi^b(T^{jm+\iota}(y'))|
  \leqslant 2^{|u|+\mathsf{p}+3} 
 \end{align*}
 and \eqref{eq: wjm long b} and \eqref{eq: wjm long c} still hold true for the same range of $\iota$. 
Finally, as an analog to \eqref{eq: wjm long final} we obtain 
\begin{align*}
\sum_{\iota=0}^{m-1}\psi^b(T^{jm+\iota}(y'))-\psi^b(T^{jm+\iota}(x'))
 &\leqslant 2^{|u|+\mathsf{p}+3}+ 6\log(2)+ 2\log(2) (q+1)m\\
 &=:U_3+2\log(2) (q+1)m. 
\end{align*}
\end{adjustwidth}

The just considered cases cover all possibilities for $j<n-1$ and finally, we consider $j=n-1$. Here, we need to consider the final entries in the definition of $x'$ in \eqref{eq: def x'}. 
Again we distinguish between different cases: 

$w_{[(n-1)m+1, nm]}=w'_{[(n-1)m+1, nm]}$ 

\begin{adjustwidth}{.8cm}{0cm}
 In all cases ($b$ being a dyadic point or not) $x'$ is constructed in such a way that 
 $$\rho(T^{(n-1)m+\iota}(x'),b)\geqslant 2^{-2m},$$ for all $\iota\in\{0,\ldots, m-1\}$. Using Lemma \ref{LEM:psi-distance-bound} implies 
 \begin{align}
  \sum_{\iota=0}^{m-1}\psi^b(T^{jm+\iota}(y'))-\psi^b(T^{jm+\iota}(x'))
  &\leqslant \sum_{\iota=0}^{m-1} 2^{2m+1-nm+(n-1)m+\iota}\leqslant 2^{2m+3}.\label{eq: last part of sum0}
 \end{align}
 Here, we applied the following consideration: Either $T^{jm+\iota}(y')$ has a larger distance to $b$ than $T^{jm+\iota}(x')$ in which case the difference is negative. Otherwise, the estimate follows from Lemma \ref{LEM:psi-distance-bound}.
 
\end{adjustwidth}

$w_{[(n-1)m+1, nm]}\in\{0^m, 1^m\}$ (in case $b$ is a dyadic point, otherwise $w_{[jm+1, (j+1)m]}=v^{m/\mathsf{p}}$ with $v\in\mathcal{V}$)
\begin{adjustwidth}{.8cm}{0cm}

This case can be dealt with more or less analogously as $w_{[jm+1, (j+1)m]}\in\{0^m, 1^m\}$ and we obtain 
\begin{align*}
 \sum_{\iota=0}^{m-1} \psi^b(T^{jm+\iota}(y'))-\psi^b(T^{jm+\iota}(x')) &\leqslant U_2-2m\leqslant U_2. 
\end{align*}
 
\end{adjustwidth}
Combining the last line with \eqref{eq: last part of sum0} yields
 \begin{align}
  \sum_{\iota=0}^{m-1}\psi^b(T^{jm+\iota}(y'))-\psi^b(T^{jm+\iota}(x'))
  &\leqslant 2^{2m+3},\label{eq: last part of sum}
 \end{align}
for $m$ sufficiently large, independent of how $w_{[(n-1)m+1, nm]}$ looks like.  
\newline

Finally, we can put the cases together. Let 
$$B_{n,m}=\{j\in \{0,\ldots, n-1\}: \mc H(w_{[jm+1, (j+q+1)m]})\cap \N_{\geqslant qm+1}\neq\emptyset\text{ for some }q\geqslant 1\}$$
and for $j\in B_{n,m}$ let $q_j=\max\{q\in \mathbb{N}:\mc H(w_{[jm+1, (j+q+1)m]})\cap \N_{\geqslant qm+1}\neq\emptyset \text{ and }j+q\leqslant n\}$. For the following, we denote the elements of $B_{n,m}$ by $j_1\leqslant\ldots\leqslant j_s$. 
Note that $j_{i+1}-j_i\geqslant q_{i}$.

We first deal with $j<j_1$. In this case we are either in case 1) or 2). Hence, 
\begin{align*}
 \sum_{\iota=0}^{j_1m-1}\psi^b(T^{\iota}(y'))-\psi^b(T^{\iota}(x'))
 &\leqslant (j_1-1)\max\{U_1,U_2\}, 
\end{align*}
for $m$ sufficiently large.
The same argument holds for $j=j+q_i+1,\ldots, j_{i+1}-1$ with $i=1,\ldots, s-1$ and we have 
\begin{align*}
 \sum_{\iota=(j_i+q_i+1)m}^{j_{i+1}m-1}\psi^b(T^{\iota}(y'))-\psi^b(T^{\iota}(x'))
 &\leqslant (j_{i+1}-(j_i+q_i)) U_1,
\end{align*}
for $m$ sufficiently large. 

Next, we deal with $j=j_i,\ldots,j_i+q_{i}$ and $i=1,\ldots, s-1$. In this case, for $j=j_i$, we are in case 3) and for $j=j_i+1,\ldots,j_i+q_i$ we are in case 2). Thus, we have 
\begin{align*}
 \MoveEqLeft\sum_{\iota=j_im}^{(j_i+q_i+1)m-1}\psi^b(T^{\iota}(y'))-\psi^b(T^{\iota}(x'))\\
 &\leqslant U_3+ 2\log(2)(q_i+1)m-2\log(2)(q_i-1)m+ (q_i-1) U_2\\
 &=U_3+ 4\log(2)m+ (q_i-1) U_2.
\end{align*}

Finally, we consider the cases $j>j_s+q_s$ if $j_s+q_s<n$. Then, using \eqref{eq: last part of sum} we have 
\begin{align*}
 \sum_{\iota=(j_s+q_i+1)m}^{j_{i+1}m-1}\psi^b(T^{\iota}(y'))-\psi^b(T^{\iota}(x'))
 &\leqslant (n-(j_s+q_s-1)) U_1+2^{3m+2}, 
\end{align*}
for $m$ sufficiently large.

Eventually, we need to estimate the whole sum where we take the maximum over all possible values for $j_1,\ldots, j_s$ and $q_1,\ldots, q_s$. There exists $U'>0$ (independent of $n$ and $m$) such that 
\begin{align*}
 \sum_{\iota=0}^{nm-1}\psi^b(T^{\iota}(y'))-\psi^b(T^{\iota}(x'))
 &\leqslant nU'+2^{3m+2}+4\log(2)m
 =nm \left(\frac{U'}{m}+\frac{2^{3m+2}+4\log(2)m}{nm}\right),
\end{align*}
for all $m$ sufficiently large. 
Choosing then first $m$ sufficiently large such that $U'/m<\varepsilon/2$ and then $n$ sufficiently large such that $(2^{3m+2}+4\log(2)m)/(nm)<\varepsilon/2$ gives the statement of the proposition. 
\end{proof}

\begin{proof}[Proof of Lemma~\ref{lem: split up of sums}]
Let $i$ be fixed and $x$ and $y$ be the points on which the supremum and infimum respectively is attained. Then Lemma \ref{LEM:psi-distance-bound} implies 
\begin{align}
 \psi^b(T^{i-1}(x))- \psi^b(T^{i-1}(y))
 &\leqslant 2^{-nm+i+1}\left(\inf_{z\in \langle\theta_m(\omega)\rangle} \rho(b,T^{i-1} (z))\right)^{-1}.\label{eq: sup psi1}
\end{align}
As $\theta_m(\omega)\in \widetilde{\Sigma}^{nm}_{m}$ we have that 
\begin{equation}
 \inf_{z\in \langle\theta_m(\omega)\rangle} \rho(b, T^{i-1}(z))\geqslant 2^{-2m+(i\mod m)-1}\geqslant 2^{-2m-1},\label{eq: sup psi2}
\end{equation}
 for all $\omega$ and $i\in \{1,\ldots, (n-1)m\}$. 
Finally, for $i\in \{nm+1,\ldots, nm\}$ using again Lemma \ref{LEM:psi-distance-bound} and the fact that $y\in \widetilde{\X}_m^b$ we have 
\begin{align}
 \psi^b(T^{i-1}(x))- \psi^b(T^{i-1}(y))
 &\leqslant - \psi^b(T^{i-1}(y))
 \leqslant -2 \log(2\rho (T^{i-1}(y), b))
 \leqslant -2 \log(2^{-2m})\notag\\
 &= 4m\log(2).
 \label{eq: sup psi3}
\end{align}
Combining \eqref{eq: sup psi1}, \eqref{eq: sup psi2} and \eqref{eq: sup psi3} implies
\begin{align*}
\MoveEqLeft\sum_{\omega\in \Sigma^{nm}} \sum_{i=1}^{nm} 
    \sup_{x\in \langle\theta_m(\omega)\rangle}\psi^b(T^{i-1}(x))- \inf_{y\in \langle\theta_m(\omega)\rangle\cap\widetilde{\X}_m^b}\psi^b(T^{i-1}(y))\\
&\leqslant \sum_{\omega\in \Sigma^{nm}} \left(\sum_{i=1}^{(n-1)m} 2^{-nm+i+1} 2^{2m+1}\right)+ \left(\sum_{i=(n-1)m+1}^{nm}4m\log(2)\right)\\
&\leqslant2^{nm}(2^{m+3}+4m^2\log(2))
\leqslant \varepsilon nm2^{nm},
\end{align*}
for $m$ fixed and $n$ sufficiently large.
\end{proof}

\begin{proof}[Proof of Lemma \ref{lem: Enm estim}]
We have
\begin{align}
 E(n,m)
 &=\frac{1}{nm}\sum_{\omega\in \Sigma^{nm}} \left(\int_{\langle\omega\rangle}S_{nm}(\overline{\psi}^b)\,\mathrm{d}\mu
 -
 \int_{\langle\theta_m(\omega)\rangle}S_{nm}(\overline{\psi}^b)\,\mathrm{d}\mu\right)\notag\\
 &=\frac{1}{nm}\sum_{\omega\in \Sigma^{nm}}\sum_{i=1}^n \left(\int_{\langle\omega_{[1,im]}\theta_m(\omega)_{[im+1,mn]}\rangle}S_{nm}(\overline{\psi}^b)\,\mathrm{d}\mu\right.\notag\\
 &\qquad\qquad\qquad\qquad\left.-
 \int_{\langle\omega_{[1,(i-1)m]}\theta_m(\omega)_{[(i-1)m+1,mn]}\rangle}S_{nm}(\overline{\psi}^b)\,\mathrm{d}\mu\right)\notag\\
 &=\frac{1}{nm}\sum_{i=1}^n\sum_{\omega\in P(i,n)} \left(\int_{\langle\omega_{[1,im]}\theta_m(\omega)_{[im+1,mn]}\rangle}S_{nm}(\overline{\psi}^b)\,\mathrm{d}\mu\right.\notag\\
 &\qquad\qquad\qquad\qquad\left.
 -
 \int_{\langle\omega_{[1,(i-1)m]}\theta_m(\omega)_{[(i-1)m+1,mn]}\rangle}S_{nm}(\overline{\psi}^b)\,\mathrm{d}\mu\right)\label{eq: Enm} 
\end{align}
where
\begin{align*}
 P(i,n)=\big\{\omega\in \Sigma^{nm}: \omega_{[im+1,(i+1)m]}\in \{u01^{m-|u|-1},u10^{m-|u|-1}\}\big\},
\end{align*}
in case $b$ is a dyadic point and 
\begin{align*}
 P(i,n)=\big\{\omega\in \Sigma^{nm}: \omega_{[im+1,(i+1)m]}=0^juv^{m/\mathsf{p}-j-|u|}, \text{ for some }j\in \{0,\ldots, \mathsf{p}-1\}\big\}
\end{align*}
otherwise, where we consider $\omega_{[i,j]}$ as the empty word if $j<i$.

Here the second equality in \eqref{eq: Enm} follows by adding summands that cancel out in a telescopic sum. The third equality follows by switching the order of summation and then leaving summands out that equal zero.

For the following let $q=u+1$ in case $b$ is a dyadic point and $q=u+\mathsf{p}$ otherwise.
We can estimate 
\begin{align}
 \MoveEqLeft\sum_{\omega\in P(i,n)} \left(\int_{\langle\omega_{[1,im]}\theta_m(\omega)_{[im+1,mn]}\rangle}S_{nm}(\psi^b)\,\mathrm{d}\mu
 -
 \int_{\langle\omega_{[1,(i-1)m]}\theta_m(\omega)_{[(i-1)m+1,mn]}\rangle}S_{nm}(\psi^b)\,\mathrm{d}\mu\right)\notag\\
 &= \sum_{\omega\in P(i,n)} \left(\int_{\langle\omega_{[1,im]}\theta_m(\omega)_{[im+1,mn]}\rangle}S_{(i-1)m+q}(\psi^b)\,\mathrm{d}\mu\right.\notag\\
 &\qquad \left.-
 \int_{\langle\omega_{[1,(i-1)m]}\theta_m(\omega)_{[(i-1)m+1,mn]}\rangle}S_{(i-1)m+q}(\psi^b)\,\mathrm{d}\mu\right)\notag\\
 &\leqslant \sum_{\omega\in P(i,n)} \left(\sup_{x\in \langle\omega_{[1,im]}\theta_m(\omega)_{[im+1,mn]}\rangle}S_{(i-1)m+q}(\psi^b)(x)\right.\notag\\
 &\qquad \left.-
 \inf_{y\in \langle\omega_{[1,(i-1)m]}\theta_m(\omega)_{[(i-1)m+1,mn]}\rangle}S_{(i-1)m+q}(\psi^b)(y)\right)2^{-nm},\label{eq: Enm estim}
\end{align}
where $2^{-nm}=\mu( \langle\omega_{[1,(i-1)m]}\theta_m(\omega)_{[(i-1)m+1,mn]}\rangle)=\mu(\langle\omega_{[1,im]}\theta_m(\omega)_{[im+1,mn]}\rangle)$.
Next, we split the sums 
\begin{align}
 \MoveEqLeft\sup_{x\in \langle\omega_{[1,im]}\theta_m(\omega)_{[im+1,mn]}\rangle}S_{(i-1)m+q}(\psi^b)(x)
 -
 \inf_{y\in \langle\omega_{[1,(i-1)m]}\theta_m(\omega)_{[(i-1)m+1,mn]}\rangle}S_{(i-1)m+q}(\psi^b)(y)\notag\\
 &\leqslant  \sum_{k=1}^{(i-1)m} \left(\sup_{x\in \langle\omega_{[k,im]}\theta_m(\omega)_{[im+1,mn]}\rangle}\psi^b(x)
 -
 \inf_{y\in \langle\omega_{[k,(i-1)m]}\theta_m(\omega)_{[(i-1)m+1,mn]}\rangle}\psi^b(y)\right)\label{eq: split of sums1}\\
 &\qquad +\sup_{x\in \langle\omega_{[(i-1)m+1,im]}\theta_m(\omega)_{[im+1,mn]}\rangle}S_q(\psi^b)(x)
 -
 \inf_{y\in \langle\theta_m(\omega)_{[(i-1)m+1,mn]}\rangle}S_q(\psi^b)(y).\label{eq: split of sums2}
\end{align}
First, we will analyze the summands of the sum \eqref{eq: split of sums1}. For this, we define the following sets
\begin{align*}
 A(i, \omega)&:= \big\{k\in \{1, \ldots, (i-1)m\}: \omega_{[k, k+m+1]}\notin \{b_{[1,m+1]}, b^{\ast}_{[1,m+1]}\}\big\}\\
  B(i,\omega)&:= \big\{k\in \{1, \ldots, (i-1)m\}: \omega_{[k, im]}\in \{b_{[1,im-k+1]}, b^{\ast}_{[1,im-k+1]}\}\big\}\\
  C(i,\omega)&:=\big\{k\in \{1, \ldots, (i-1)m\}: \omega_{[k, k+m+1]}\in \{ b_{[1,m+1]}, b^{\ast}_{[1,m+1]}\}\\
  &\hphantom{=C(i,\omega)=\{k\in \{1, \ldots, (i-1)m\}:} \text{ but }\omega_{[k, im]}\notin \{b_{[1,im-k+1]}, b^{\ast}_{[1,im-k+1]}\}\big\}.
\end{align*}
It is clear that $A(i, \omega)\cup B(i, \omega)\cup C(i, \omega)=\{1,\ldots, (i-1)m\}$. 
Let us analyze the sums separately. For each $k\in A(i,\omega)$ we have by Lemma \ref{LEM:psi-distance-bound}
\begin{align*}
 \sup_{x\in \langle\omega_{[k,im]}\theta_m(\omega)_{[im+1,mn]}\rangle}\psi^b(x)
 -
 \inf_{y\in \langle\omega_{[k,(i-1)m]}\theta_m(\omega)_{[(i-1)m+1,mn]}\rangle}\psi^b(y)
 &\leqslant \frac{2}{2^{-m-1}}2^{k-(i-1)m+1}\\
 &=2^{k-im+3}.
\end{align*}
Summing over all $k\in A(i,\omega)$ we obtain 
\begin{align}
 \MoveEqLeft\sum_{k\in A(i,\omega)} \left(\sup_{x\in \langle\omega_{[k,im]}\theta_m(\omega)_{[im+1,mn]}\rangle}\psi^b(x)
 -
 \inf_{y\in \langle\omega_{[k,(i-1)m]}\theta_m(\omega)_{[(i-1)m+1,mn]}\rangle}\psi^b(y)\right)\notag\\
 &\leqslant \sum_{k=1}^{(i-1)m} 2^{k-(i-2)m+3}=2^{m+3}.\label{eq: sum estim A} 
\end{align}
Next, we analyze $k\in B(i,\omega)$. It is clear that $\#B(i,\omega)\in \{0,1\}$. Moreover, we have for $k\in B(i,\omega)$ that
\begin{align*}
 \sup_{x\in \langle\omega_{[k,im]}\theta_m(\omega)_{[im+1,mn]}\rangle}\psi^b(x)
 <
 \inf_{y\in \langle\omega_{[k,(i-1)m]}\theta_m(\omega)_{[(i-1)m+1,mn]}\rangle}\psi^b(y)
\end{align*}
and thus 
\begin{align}
 \MoveEqLeft\sum_{k\in B(i,\omega)} \left(\sup_{x\in \langle\omega_{[k,im]}\theta_m(\omega)_{[im+1,mn]}\rangle}\psi^b(x)
 -
 \inf_{y\in \langle\omega_{[k,(i-1)m]}\theta_m(\omega)_{[(i-1)m+1,mn]}\rangle}\psi^b(y)\right)< 0. \label{eq: sum estim B}
\end{align}
Finally, we analyze $k\in C(i,\omega)$.
Let $h(k)=\max\{\ell\in \mathbb{N}: \omega_{[k,k+\ell-1]}\in \{b_{[1,\ell]}, b^{\ast}_{[1,\ell]}\}\}$. 
By Lemma \ref{LEM:psi-distance-bound} we have for $k\in C(i,\omega)$ that  
\begin{align*}
 \sup_{x\in \langle\omega_{[k,im]}\theta_m(\omega)_{[im+1,mn]}\rangle}\psi^b(x)
 -
 \inf_{y\in \langle\omega_{[k,(i-1)m]}\theta_m(\omega)_{[(i-1)m+1,mn]}\rangle}\psi^b(y)
 &\leqslant 2^{-(i-1)m+k+1+h(k)+1}.
\end{align*}
Let $\rho_1< \rho_2<\ldots < \rho_\ell$ be all elements of $C(i,\omega)$. Then we have $h(\rho_j)\leqslant \rho_{j+1}-\rho_j$, for all $j\in \{1,\ldots, \ell-1\}$. 

Moreover, $\omega_{[k, im]}\notin \{b_{[1,im-k+1]}, b^{\ast}_{[1,im-k+1]}\}$ for all $k\in \{\rho_1,\ldots, \rho_{\ell}\}$
immediately implies that 
$\omega_{[\rho_{\ell}, i(m-1)+1]}\notin \{b_{[1,i(m-1)-k+2]}, b^{\ast}_{[1,i(m-1)-k+2]}\}$, as we only consider $\omega_{[i(m-1)+1, im]}\in \{0^m, 1^m\}$ in the dyadic case and $\omega_{[i(m-1)+1, im]}= v^{m/\mathsf{p}}$ for $v\in\mathcal{V}$ otherwise. This in particular implies 
$h(\rho_{\ell})\leqslant (m-1)i-\rho_{\ell}$. 
Hence, 
\begin{align}
  \MoveEqLeft\sum_{k\in C(i,\omega)} \left(\sup_{x\in \langle\omega_{[k,im]}\theta_m(\omega)_{[im+1,mn]}\rangle}\psi^b(x)
 -
 \inf_{y\in \langle\omega_{[k,(i-1)m]}\theta_m(\omega)_{[(i-1)m+1,mn]}\rangle}\psi^b(y)\right)\notag\\
 &\leqslant \left(\sum_{j=1}^{\ell-1}2^{-(i-1)m+\rho_{j}+1+\rho_{j+1}-\rho_{j}}\right)+2^{-(i-1)m+\rho_{\ell}+1+(m-1)i-\rho_{\ell}}\notag\\
 &=\left(\sum_{j=1}^{\ell-1}2^{-(i-1)m+\rho_{j+1}+1}\right)+2^m\leqslant 8+2^m\leqslant 2^{m+1}, \label{eq: sum estim C}
\end{align}
for $m$ sufficiently large. 

Combining \eqref{eq: sum estim A}, \eqref{eq: sum estim B}, and \eqref{eq: sum estim C}
implies 
\begin{align}
 \MoveEqLeft\sum_{k=1}^{(i-1)m} \left(\sup_{x\in \langle\omega_{[k,im]}\theta_m(\omega)_{[im+1,mn]}\rangle}\psi^b(x)
 -
 \inf_{y\in \langle\omega_{[k,(i-1)m]}\theta_m(\omega)_{[(i-1)m+1,mn]}\rangle}\psi^b(y)\right)
 \leqslant 2^{m+4},\label{eq: split of sums1estim}
\end{align}
for $m$ sufficiently large.

The second sum \eqref{eq: split of sums2} can be estimated as 
\begin{align}
 \MoveEqLeft\sup_{x\in \langle\omega_{[(i-1)m+1,im]}\theta_m(\omega)_{[im+1,mn]}\rangle}S_q(\psi^b)(x)
 -
 \inf_{y\in \langle\theta_m(\omega)_{[(i-1)m+1,mn]}\rangle}S_q(\psi^b)(y)\notag\\
 &\leqslant 0-\inf_{y\in \mathcal{D}}S_q(\psi^b)(y)=: C, \label{eq: split of sums2estim}
\end{align}
where 
\begin{align*}
 \mathcal{D}=\begin{cases}\langle 0^{q}\rangle \cup \langle 1^{q}\rangle&\text{if }b\text{ is a dyadic point}\\
    \bigcup_{v\in \mathcal{V}}v^{q/\mathsf{p}}&\text{otherwise}
   \end{cases}
\end{align*}
and $C$ is independent of $m$ and $n$. 

Inserting the estimates \eqref{eq: split of sums1estim} and \eqref{eq: split of sums2estim} in \eqref{eq: split of sums1}/\eqref{eq: split of sums2}
and combining this with \eqref{eq: Enm} and \eqref{eq: Enm estim}
yields 
\begin{align*}
 E(n,m)
 &\leqslant \frac{1}{nm}\sum_{i=1}^n\sum_{\omega\in P(i,n)} (2^{m+4}+C)  2^{-nm},
\end{align*}
for $m$ sufficiently large. Noting that $\#P(i,n)=2^{m(n-1)+1}$ and $2^{m+4}+C<2^{m+5}$, for $m$ sufficiently large, implies 
$E(n,m)\leqslant 2^6/m\leqslant \varepsilon$, for $m$ sufficiently large.
\end{proof}

\section{Proof of the large deviation results}\label{sec: proofs large dev}

\subsection{Proof of one sided large deviation (Lemma \ref{lem:upper_bound})}
\begin{proof}[Proof of Lemma \ref{lem:upper_bound}]
For $M>0$ define a truncation $\psi_M^b$
by
\[ \psi_M^b (x) = \left\{ \begin{array}{ll}
         \psi^b (x)  & \mbox{if $\psi^b (x) \geqslant -M$};\\
        0 & \mbox{if $\psi^b (x) < -M$}.\end{array} \right. 
        \]

 Choose $M$ large enough that $|\mu(\psi_M^b)-\mu (\psi^b)|<\varepsilon$, and note that $
 0<\mu(\psi_M^b)-\mu (\psi^b)<\varepsilon$. 
 The function $\psi_M^b$ is bounded and of bounded variation.
 Hence  there exists $C_0(\varepsilon)>0$ and a function $g_0 (\varepsilon)>0$ such that for all $n\geqslant 0$,
 \[
 \mu ( S_n(\psi_M^b-\mu(\psi_M^b))>n\varepsilon)<C_0e^{-ng_0(\varepsilon)}.
 \]
 Since $S_n (\psi_M^b-\mu(\psi^b))=S_n(\psi_M^b-\mu(\psi_M^b)) +n (\mu(\psi_M^b)-\mu (\psi^b))$ we have 
 $$\left\{ S_n(\psi_M^b-\mu(\psi^b))>2n\varepsilon\right\}  \subset \left\{ S_n(\psi_M^b-\mu(\psi_M^b))>n\varepsilon\right\}$$ and so
 \[
 \mu ( S_n(\psi_M^b-\mu(\psi^b))>2n\varepsilon)<C_0e^{-ng_0(\varepsilon)}.
 \]
 Since $S_n (\psi_M^b)\geqslant S_n (\psi^b)$ we see 
 \[
 \mu ( S_n(\psi^b-\mu(\psi^b))>2n\varepsilon)<Ce^{-ng_0(\varepsilon)}
 \]
 which implies the result by taking $g(\varepsilon)=g_0(\varepsilon/2)$ and $C(\varepsilon)=C_0(\frac{\varepsilon}{2})$.
 \end{proof}

\subsection{Bounds for $S_n^I (x)$ from below via the pressure function - proof of Lemma \ref{inf:bound}}\label{inf_bound}

In the next lemma we detail the consequences of the limit in \eqref{eq: pressure diff}. Recall $\overline{\psi}^b=\psi^b-\mu(\psi^b)$.

 \begin{lemma}\label{pressure_limit}
 If $\lim_{t\uparrow 0}\frac{1}{t}[\mathcal{P}_{\operatorname{top}}(t\overline{\psi}^b)-\mathcal{P}_{\operatorname{top}} (0)]=0$, then  given $\varepsilon>0$,
$|\mathcal{P}_{\operatorname{top}} (-q\overline{\psi}^b)  -\log (2)|< q\varepsilon$ for  all small enough $-q<0$.
\end{lemma}
 
\begin{proof}
Assume that  $\lim_{t\uparrow 0}[\mathcal{P}_{\operatorname{top}}(t\overline{\psi}^b)-\log (2)]=o(t)$. Given $\eta>0$ there
exists $\delta_{\eta}>0$ such that $|\mathcal{P}_{\operatorname{top}}(t\overline{\psi}^b) -\log (2)| < |\eta t|$ for $|t|<\delta_{\eta}$, $t<0$.

Given $\varepsilon >0$ choose $0<\eta <\varepsilon$.
Then   $|\mathcal{P}_{\operatorname{top}}(-q\overline{\psi}^b) -\log (2)|\leqslant |q\eta|  <|q\varepsilon|$ for small enough $|q|< \delta_{\eta}$,$-q<0$.
\end{proof}

Since $\overline{\psi}^b=\psi^b-\mu (\psi^b)$  we have  $\mu ( \overline{\psi}^b) =0$.

For the following, as a slight abuse of notation we also write $S_n^I (\overline{\psi}^b) (\langle \omega\rangle)$ when we mean $S_n^I (\overline{\psi}^b) (x)$ for any $x\in\langle \omega\rangle$.

\begin{proof}[Proof of Lemma \ref{inf:bound}] 
Let $J_n^{-} (\varepsilon) =\{ \omega\in \Sigma^n : S_n^I (\overline{\psi}^b) (x) \leqslant -n\varepsilon \text{ for }x\in \langle \omega\rangle\}$. 
Using the inequality $-qS^I_n (\overline{\psi}^b(\langle \omega\rangle)+\varepsilon) \geqslant 0$ if $\omega\in J_n^-(\varepsilon)$, we obtain for $q>0$ 
\begin{align*}
 \mu \Big(\bigcup_{\omega\in J_n^{-}(\varepsilon) }\langle \omega\rangle\Big)
&\leqslant \sum_{\omega\in J_n^{-}(\varepsilon) } e^{-n\log (2)-q S_n^I (\overline{\psi}^b (\langle \omega\rangle)+n\varepsilon)}\\
&= \sum_{\omega\in J_n^{-}(\varepsilon)} 2^{-n}e^{-nq\varepsilon}e^{ -qS_n^I(\overline{\psi}^b)(\langle \omega\rangle)}\\
&\leqslant  2^{-n}e^{-nq\varepsilon} \sum_{\omega\in \Sigma^n}e^{-qS_n^I(\overline{\psi}^b)(\langle \omega\rangle)}.
\end{align*}
Thus,
\[
 \frac{1}{n}\log \Big(\mu \Big(\bigcup_{\omega\in J_n^{-}(\varepsilon) }\langle \omega\rangle\Big)\Big)\leqslant -\log (2) -q\varepsilon +\frac{1}{n} \log \Big(\sum_{\omega \in \Sigma^n} e^{-qS_n^I(\overline{\psi}^b)(\langle \omega\rangle)}\Big).
\]
Hence, 
\[
\limsup_{n\to \infty} \frac{1}{n}\log \Big(\mu \Big(\bigcup_{\omega\in J_n^{-}(\varepsilon) }\langle \omega\rangle\Big)\Big)\leqslant -\log  (2) - q\varepsilon + \mathcal{P}_{\operatorname{top}}(-q(\overline{\psi}^b)).
\]
Since, by Lemma \ref{pressure_limit}, we have $|\mathcal{P}_{\operatorname{top}}(- q(\overline{\psi}^b))-\log (2)| <q\frac{\varepsilon}{2} $ for all sufficiently small $q>0$, this implies the result.
\end{proof}

\subsection{Exponential  bounds for $S_n(\overline{\psi}^b)$ - proof of Lemma~\ref{lem:cylinder-to-orbit}}\label{exp_bounds}

Before starting with the proof of Lemma \ref{lem:cylinder-to-orbit} we will first need to prove the following lemma:

\begin{lemma}\label{lemma_upper_bound}
Given $\varepsilon>0$,
\begin{equation}\label{ep:upper_bound}
\mu \left( x: S_n^I (\overline{\psi}^b)(x) \geqslant n\varepsilon \right) \leqslant C_0(\varepsilon)e^{-ng_0(\varepsilon)}
\end{equation}
 for some $C_0(\varepsilon)>0$, $g_0 (\varepsilon)>0$.
\end{lemma}

\begin{proof}
Let $\varepsilon>0$ and choose $0<\beta<1$ such that $|\beta\mu(\psi^b)|<\frac{\varepsilon}{4}$.
For each $x$, let $q_n(x)$ be a point in $\langle x_{[1,n]}\rangle$, such that $S_n^I (\overline{\psi}^b) (x)=
\sum_{j=0}^{n-1} (\overline{\psi}^b\circ T^j)  (q_n(x))$.

Let $B_n:=\{x: |x-b|\leqslant 2^{-\frac{n\beta}{2}}\}$.  By the mean value theorem, as $|D\psi^b (x)| \leqslant \frac{C}{|x-b|}$ for some $C>0$ and  $|x-b|$ small we have that if 
$T^j (\langle x_{[1,n]}\rangle) \cap B_n=\emptyset$ for $j=0,1,...,\lfloor (1-\beta)n\rfloor$, then
\[
|\psi^b (T^j (x))-\psi^b  (T^j( q_n(x)))| \leqslant C2^{\frac{n\beta}{2}}   2^{-(n-j)},
\]
for $j=0,1,...,\lfloor (1-\beta)n\rfloor$.
Hence,
\[
\sum_{j=0}^{\lfloor (1-\beta)n\rfloor} \left|\psi^b (T^j (x))-\psi^b  (T^j (q_n(x)))\right|\leqslant C2^{\frac{n\beta}{2}} 2^{-n\beta+1}
=C2^{-\frac{n\beta}{2}+1}<\frac{n\varepsilon}{4},
\]
for $x\in B_n$ and $n$ sufficiently large.

 We now estimate the measure of cylinder sets $\langle \omega\rangle$,  $\omega\in \Sigma^n$ with the property that $T^{j}(\langle \omega\rangle) \cap B_n \not = \emptyset$
 for any $j=1,\ldots,\lfloor (1-\beta)n\rfloor$. For each $j=1,...,\lfloor (1-\beta)n\rfloor$, we have that $\{T^j(\langle \omega\rangle): \omega\in \Sigma^n\}$ is a partition
 of $[0,1]$ into intervals of length $2^{-(n-j)}$, where $2^{-(n-j)} \leqslant 2^{-n(1-\beta)}$. 
 The number  of cylinders 
 in $\{T^j(\langle \omega\rangle): \omega\in \Sigma^n\}$ intersecting $B_n$ is at most $\lceil 2^{-\frac{n\beta}{2}} 2^{n-j}\rceil=\lceil 2^{n(1-\frac{\beta}{2} )-j}\rceil$, 
 each of measure
 $2^{-(n-j)}$, so the total fraction of $\omega \in \Sigma^n$, such that $T^j (\langle \omega\rangle)$  intersects $B_n$ is 
 $\lceil 2^{n(1-\frac{\beta}{2} )-j}\rceil 2^{-(n-j)}\leqslant 2^{-\frac{n\beta}{2}}+2^{-(n-j)}$. Hence,
 the measure of all cylinder sets $\langle \omega\rangle$, $\omega\in \Sigma^n$ with the property that $T^{j}(\langle \omega\rangle) \cap B_n \not = \emptyset$
 for any $j=1,\ldots,\lfloor(1-\beta)n\rfloor$ is at most $(1-\beta)n2^{-\frac{n\beta}{2}}$, for $n$ sufficiently large.

Thus, except for a set of measure  at most $(1-\beta)n2^{-\frac{n\beta}{2}}$ we have 
$$\sum_{j=0}^{\lfloor (1-\beta)n\rfloor} \left|\psi^b (T^j (x))-\psi^b  (T^j (q_n(x)))\right|< \frac{n\varepsilon}{4}.$$

Moreover, by Lemma~\ref{lem:upper_bound} applied with $\frac{\varepsilon}{4}$ we have
\[
\sum_{j=0}^{\lfloor(1-\beta)n\rfloor} \left(\psi^b (T^j(x) ) -\mu(\psi^b) \right)\leqslant \frac{n\varepsilon}{4}
\]
except for a set of measure at most $C(\varepsilon/4) e^{-(1-\beta)nI (\varepsilon/4)}$ yielding
\[
\sum_{j=0}^{\lfloor(1-\beta)n\rfloor} \left(\psi^b (T^j(q_n(x)) ) -\mu(\psi^b) \right)\leqslant \frac{n\varepsilon}{2}
\]
except for a set of measure at most $C(\varepsilon/4) e^{-(1-\beta)nI (\varepsilon/4)}+(1-\beta)n2^{-\frac{n\beta}{2}}$.

As a simple upper bound, we obtain, using $\beta |\mu(\psi^b)|< \frac{\varepsilon}{4}$ that 
\[
\sum_{j=\lfloor (1-\beta)n\rfloor+1}^n \psi^b  (T^j (q_n(x))) -\mu(\psi^b) \leqslant (n-\lfloor (1-\beta)n\rfloor) |\mu(\psi^b)|\leqslant n\frac{\varepsilon}{2},
\]
for $n$ sufficiently large. 
So, for $n$ sufficiently large, we have 
\[
\sum_{j=0}^n \left(\psi^b (T^j (q_n(x))) - \mu(\psi^b)\right) 
\leqslant n\varepsilon
\]
except for a set of measure less than $C(\varepsilon/4) e^{-(1-\beta)nI (\varepsilon/4)}+(1-\beta)n2^{-\frac{n\beta}{2}}$. 
Thus
\[
\mu \left( x: S_n^I (\overline{\psi}^b) (x)  \geqslant n\varepsilon\right)\leqslant C_0(\varepsilon) e^{-ng_0(\varepsilon)}
\]
for some $g_0(\varepsilon)>0$, $C_0(\varepsilon)>0$.
\end{proof}

\begin{proof}[Proof of Lemma \ref{lem:cylinder-to-orbit}]

For each $x$, let $q_n(x)$ be a point in $\langle x_{[1,2n]}\rangle$, such that $S_{2n}^I (\overline{\psi}^b) (x)=
\sum_{j=0}^{2n-1} (\overline{\psi}^b\circ T^j)  (q_n(x))$.
Our first claim is that there exist  $C_3(\varepsilon)$, $g_3(\varepsilon)$ such that
\begin{equation}\label{eq: half sum > 3neps}
 \mu\left(\sum_{j=0}^{n -1} (\overline{\psi}^b \circ T^j) ( q_n(x)) \leqslant -6 n\varepsilon\right)\leqslant C_3(\varepsilon)e^{-ng_3(\varepsilon)}.
\end{equation}
To see this we note that by definition 
\[
S_{2n}^{I} (\overline{\psi}^b) (x)=\sum_{j=0}^{2n-1} (\overline{\psi}^b \circ T^j) ( q_n(x)) =\sum_{j=0}^{n -1} (\overline{\psi}^b \circ T^j)  (q_n(x))
+ \sum_{j=n}^{2n-1} (\overline{\psi}^b \circ T^j)  (q_n(x)).
\]
For the second summand we have
\[
\sum_{j=n}^{2n-1} (\overline{\psi}^b\circ T^j)  (q_n(x)) 
\leqslant S_{n}^I (\overline{\psi}^b) (T^{n} (x))
\]
and by Lemma \ref{lemma_upper_bound} we have that $S_{n}^I (\overline{\psi}^b) (T^{n}  (x)) \leqslant \varepsilon n$ except for a set of measure
at most $C_0 (\varepsilon)e^{-ng_0 ( \varepsilon)}$.

Hence, if $ \sum_{j=0}^{n -1} (\overline{\psi}^b \circ T^j)  (q_n(x)) \leqslant -6 n\varepsilon$ then except for a set of measure at most $C_0 (\varepsilon)e^{-ng_0 ( \varepsilon)}$ we have
\[
S_{2n}^I (\overline{\psi}^b) (x) < -6n\varepsilon+ n\varepsilon <-4 n\varepsilon
\]
which means that $x\in \{S_{2n}^I (\overline{\psi}^b)  \leqslant -2n \varepsilon\}$.
However, under the assumptions of the lemma we have
\begin{align*}
 \mu \left(x:S_{2n}^I (\overline{\psi}^b) (x)  \leqslant -2n \varepsilon\right)
 &\leqslant C_1 (\varepsilon) e^{-2n g_1(\varepsilon)}.\label{eq: SnI<-neps}
 \end{align*}
Thus except for a set of measure less than $C_3 (\varepsilon)e^{-n g_3 (\varepsilon)}$ we have 
$\sum_{j=0}^{n -1} (\overline{\psi}^b \circ T^j)  (q_n (x)) \geqslant -6 n\varepsilon$ where $C_3 (\varepsilon)=C_0 (\varepsilon) +C_1 (\varepsilon) >0$ and $g_3 (\varepsilon)=\min \{g_0 ( \varepsilon),2g_1 (\varepsilon)\}>0$ establishing \eqref{eq: half sum > 3neps}.

Consider the set 
   \[
  B_n:=\{x:|x-b| \leqslant 2^{-n/4}\}.
  \]
  By the mean value theorem, there is a constant $C>0$ such that  $|D\overline{\psi}^b (x)| \leqslant \frac{C}{|x-b|}$ for  $|x-b|$ small. Hence if 
  $T^j \langle x_{[1,2n]}\rangle \cap B_{n}=\emptyset$ for $j=0,1,...,n -1$ then
\[
|\overline{\psi}^b (T^j (x))- \overline{\psi}^b  (T^j (q_n(x)))|\leqslant 2^{n/4}2^{-(2n-j)}
\]
for $j=0,1,...,n-1$.

 We now estimate the measure of cylinder sets  $\langle w\rangle$, $w\in \{0,1\}^{2n}$ with the property that $T^{j}(\langle w\rangle) \cap B_n \not = \emptyset$
 for any $j=0,\ldots,n-1$. 
 Since for each $j=0,...,n-1$, $\{T^j (\langle w\rangle) : w \in \{0,1\}^n\}$ is a partition
 of $[0,1]$ into intervals of length $2^{-(2n-j)} < 2^{-n/4}$, so as an upper estimate we obtain 
 $$\mu\left( T^j (\langle w\rangle): w\in \{0,1\}^{2n}, T^j \langle w\rangle\cap B_n \not =\emptyset
 \mbox{ for any } j=0,...,n-1\right)\leqslant n 2^{-n/4}. $$

Thus except for a set of measure  at most $n 2^{-n/4} $ we have 
\[
\sum_{j=0}^{n-1 } |\overline{\psi}^b (T^j (x))-\overline{\psi}^b  (T^j (q_n(x)))|<n 2^{-n/4}.
\]
So we have 
\[
\sum_{j=0}^{n-1} \overline{\psi}^b (T^j (q_n(x))) \leqslant \sum_{j=0}^{n-1} \overline{\psi}^b (T^j(x) ) +n 2^{- n/4}
\]
except for a set of measure less than $n 2^{-n/4}$.

Thus, by \eqref{eq: half sum > 3neps}, for $n$ sufficiently large and except for a set of measure less than $n 2^{-n/4}+C_3(\varepsilon)e^{-ng_3( \varepsilon)}$ we have
\[
\{ x: S_{n}(\overline{\psi}^b)(x)  \leqslant -8n\varepsilon\}
\subset \bigg\{ x: \sum_{j=0}^{2n-1} \overline{\psi}^b(T^j (q_n(x))) \leqslant -6n \varepsilon\bigg\}
\]
and so 
\[
\mu \left( x: S_{n}(\overline{\psi}^b)(x) \leqslant -8 n\varepsilon\right)\leqslant C_4(\varepsilon) e^{-ng_4(\varepsilon)}
\]
for constants $C_4(\varepsilon)>0$ and $g_4 (\varepsilon)>0$. Since $\varepsilon$ was arbitrary we may define $C_2(\varepsilon)=C_4 (\varepsilon/8)$ and $g_2(\varepsilon)= 2g_4 (\varepsilon/8)$ to finish the proof.
\end{proof}

\subsection{Comparing the large deviation of $\phi^b$ and $\psi^b$ - Proof of Lemma \ref{lem:psi-to-phi}}

\begin{proof}[Proof of Lemma \ref{lem:psi-to-phi}]
Let us assume that $b\neq 0$ (in which case by \cite{NT_LD} we know that we don't have exponential large deviations). 
Then both functions $\phi^b$ and $\psi^b$ are bounded outside a neighborhood of $b$.
The reason we exclude $b=0$ is that the singularity of $\psi^b$ is then two-sided.
So we may define a continuous (hence bounded)  function $g_b(x)$ on $[0,1]$  by
\[ g_b(x) = \left\{ \begin{array}{ll}
     \phi^b(x)   -\psi^b (x)  & \mbox{if $x \neq b$};\\
        -2\log (\pi) & \mbox{if $x =b$}.\end{array} \right. \]

We have that 
\begin{equation}
\phi^b (x) -\mu(\phi^b)= [g_b (x) - \mu(g_b)] +[\psi^b (x)-\mu(\psi^b)].\label{eq: compare phipsi}
\end{equation}
If $x\neq b$ this immediately follows from the definition of $g_b$. On the other hand, for  fixed $b\neq 0$, we have
\begin{align*}
  \lim_{x\to b}(\phi^b(x)-\psi^b (x)) =\lim_{x\to b} 2\log \left(\frac{|x-b|}{|\sin(\pi (x-b))|}\right)=-2\log (\pi).
 \end{align*} 
From \eqref{eq: compare phipsi} we immediately obtain
\[
\mu(|S_n (\phi^b-\mu(\phi^b)) |\geqslant n\varepsilon)\leqslant \mu\Big(|S_n (g_b-\mu(g_b)) |\geqslant \frac{n\varepsilon}{2}\Big)+\mu\Big(|S_n (\psi^b-\mu(\psi^b)) |\geqslant \frac{n\varepsilon}{2}\Big).
\]
Since $g_b$ is a continuous function on the doubling map, for which Lebesgue measure is the unique
equilibrium measure for the geometric potential $-\log (DT)=-\log (2)$,  by results of Young~\cite[Theorem 1 (2)]{You90} (in the notation of that theorem, taking $c=\varepsilon$, $\zeta=\log (2)$ and taking once $g_b-\mu(g_b)$ and once $\mu(g_b)-g_b$ as the observable),
we obtain
\[
\mu\Big(|S_n (g_b-\mu(g_b)) |\geqslant \frac{n\varepsilon}{2}\Big)\leqslant C_0(\varepsilon/2) e^{-nI_0(\varepsilon/2)}
\]
where $I_0 (\varepsilon)>0$ is  defined on a neighborhood $U$ of $0$, $C_0 (\varepsilon)>0$  and $I_0(\varepsilon)>0$ for $\varepsilon >0$. $I_0(\varepsilon)$ is not necessarily a rate function (by which we mean strictly convex on a neighborhood of $0$).

Let $U$ be the neighborhood of $0$ described above for $g_b$ and $U_1$ that  for $\psi^b$ on which a large deviation result holds.  
If $\varepsilon \in U\cap U_1=:U_2$,
$\varepsilon\not =0$,  then $\varepsilon/2\in U\cap U_1$.
Hence, if $\varepsilon \in U_2$ we may take $C_2=\max \{C_0,C_1\}$ and $I_2 (\varepsilon)=\min \{I_0 (\frac{\varepsilon}{2}), I_1 (\frac{\varepsilon}{2})\}>0$.
\end{proof}

\section{Proof of the case for the singularity at a periodic point (Proposition~\ref{prop: periodic part})}\label{sec: periodic case}

\begin{proof}[Proof of Proposition \ref{prop: periodic part}]
 The lower bound follows as in \cite[Theorem 3.4]{NT_LD}. 
  For the upper bound we also use the same approach as in \cite[Theorem 3.4]{NT_LD}.   
  For $0<\beta<1$ we define a sequence of functions (truncated
  $\phi_{\mathsf{p}}$) by
  \[
    g_n(x) = \left\{ \begin{array}{ll}
                       0 & \mbox{if $ d(x,\mathsf{p}) \leqslant e^{-n^{\beta}}$};\\
                       \varphi_{\mathsf{p}} (x) & 
              \mbox{if $d(x,\mathsf{p}) > e^{-n^{\beta}}$}.\end{array}
     \right.
  \]
  
 Then define $G_n:=\sum_{j=0}^{n-1} g_n\circ T^j$.

 However, we apply the  improved version \cite[Lemma 4.5]{KessSchi} instead of \cite[Lemma 6.15]{Tanja_PhD} and obtain 
 with the notation as in \cite{NT_LD} for $n$ sufficiently large that
 \begin{align*}
  \mu(|G_n - \mu(G_n)| > \varepsilon \mu(G_n) ) \leqslant  K 
  \exp\left( -\frac{\varepsilon n \mu(g_n) }{n^{\beta\alpha}}\right)
  \leqslant C_0 \exp\left(-\varepsilon n^{1-\beta\alpha}\right),
 \end{align*}
 as $\mu(g_n) \geqslant C>0$ for $n$ sufficiently large. 
 Applying  again the estimate $\mu(S_n (\varphi_{\mathsf{p}}) \neq  G_n) \leqslant n\exp(-n^{\beta})$ and setting $\beta=1/(1+\alpha)$ yields 
 \begin{align*}
  \mu\left(S_n(\phi) - n\int \phi\, \mathrm{d}\mu \geqslant  n\varepsilon\right)
  &\leqslant K \exp\left(-\varepsilon n^{\frac{1}{1+\alpha}}\right)
  +  n\exp\left(-n^{\frac{1}{1+\alpha}}\right)
  \leqslant 2K \exp\left(-\varepsilon n^{\frac{1}{1+\alpha}}\right),
 \end{align*}
 for $n$ sufficiently large.
\end{proof}

For the purpose of open access, the authors have applied a ``Creative
Commons Attribution (CC BY) licence'' to any Author Accepted Manuscript
version arising from this submission.

\end{document}